\documentclass[11pt]{article}
\usepackage{amssymb}
\usepackage{amsmath}
\usepackage{latexsym}
\usepackage{hyperref}

\usepackage[all]{xy}
\usepackage{color}

\numberwithin{equation}{section}
\newtheorem{thm}[equation]{Theorem}
\newtheorem{prop}[equation]{Proposition}
\newtheorem{defn}[equation]{Definition}
\newtheorem{rem}[equation]{Remark}
\newtheorem{lem}[equation]{Lemma}
\newtheorem{corol}[equation]{Corollary}
\newtheorem{example}[equation]{Example}
\title{Higher order elliptic operators on variable domains. Stability results and boundary oscillations for intermediate problems.   }
\author{Jos\'{e} M. Arrieta and Pier Domenico Lamberti }

\date{\ }

\begin{document}

\newcommand{\rea}{\mathbb{R}}

\newcommand{\bN}{\mathbb{N}}
\newcommand{\bR}{\mathbb{R}}
\newcommand{\bC}{\mathbb{C}}
\newcommand{\eps}{\epsilon}
\newcommand{\veps}{\varepsilon}
\newcommand{\vfi}{\Phi}
\newcommand{\tw}{\widetilde{w}}
\newcommand{\tv}{\widetilde{v}}
\newcommand{\tQ}{\widetilde{Q}}
\newcommand{\tPhi}{\widetilde{\Phi}}
\newcommand{\tmu}{\widetilde{\mu}}
\newcommand{\tq}{\widetilde{q}}
\newcommand{\txi}{\widetilde{\xi}}
\newcommand{\tdelta}{\widetilde{\delta}}
\newcommand{\bw}{\bar{w}}
\newcommand{\bv}{\bar{v}}
\newcommand{\uv}{\underline{v}}
\newcommand{\df}{\delta f}
\newcommand{\cL}{\mathcal{L}}
\newcommand{\cLinv}{\mathcal{L}^{-1}}
\newcommand{\real}{{\rm Re}\,}
\newcommand{\imag}{{\rm Im}\,}
\newcommand{\sgn}{\,{\rm sgn}\,}
\newcommand{\cC}{\,{\mathcal C}\,}
\newcommand{\vfila}{\varphi_{\lambda}}
\newcommand{\la}{\lambda}
\newcommand{\cs}{c^*}

\newcommand{\cP}{\,{\mathcal P}\,}
\newcommand{\cQ}{\,{\mathcal Q}\,}
\newcommand{\cPQ}{\,{\mathcal PQ}\,}
\newcommand{\cE}{\,{\mathcal E}\,}

\newcommand{\eto}{\stackrel{\eps\to 0}{\longrightarrow}}

\newcommand{\weto}{\stackrel{\eps\to 0}{\rightharpoonup}}

\newcommand{\cEE}{\,{\mathcal{EE}}\,}
\newcommand{\Pcon}{\stackrel{\cP}{\rightarrow}}
\newcommand{\Qcon}{\stackrel{\cQ}{\rightarrow}}
\newcommand{\PQcon}{\stackrel{\cPQ}{\longrightarrow}}
\newcommand{\Econ}{\stackrel{\cE}{\rightarrow}}
\newcommand{\EEcon}{\stackrel{\cEE}{\rightarrow}}
\newcommand{\Ccon}{\stackrel{\cC}{\rightarrow}}

\maketitle

%
%
%

\noindent 
{\bf Abstract:}  We  study the spectral behavior of higher order elliptic operators upon  domain perturbation.  We prove general spectral stability results for Dirichlet, Neumann and intermediate boundary conditions. Moreover, we consider the case of the bi-harmonic operator with those intermediate boundary conditions which appears in 
 study of hinged plates.  In this case, we analyze the spectral behavior when the boundary of the domain is subject to a periodic oscillatory perturbation.  We will show that there is a critical oscillatory behavior and the limit problem depends on whether we are above, below or just sitting on this critical value.
In particular, in the critical case we identify the strange term appearing in the limiting boundary conditions by using the unfolding method from homogenization theory. 

\vspace{11pt}

\noindent
{\bf Keywords:}
higher order elliptic operators, Dirichlet, Neumann, intermediate boundary conditions,  oscillatory boundaries, homogenization

\vspace{6pt}
\noindent
{\bf 2000 Mathematics Subject Classification:}  35J40, 35B20, 35B27,  35P15
%

%
%

\section{Introduction}
\label{introduction}

In this paper, we consider the general problem of the spectral behavior  of  an elliptic partial differential   operator  (i.e., the behavior of its eigenvalues, eigenfunctions as well as of the solutions to the corresponding Poisson problem)  when the underlying domain is perturbed. In ${\mathbb{R}}^N$ with $N\geq 2$,  we will consider a family of domains $\{\Omega_\eps\}_{0<\eps\leq \eps_0}$ which approach a limiting domain $\Omega$ as $\eps\to 0$,  in certain sense to be specified and we will also consider higher order selfadjoint operators (order $2m$ with $m\geq 1$) with not necessarily constant coefficients and with certain boundary conditions.  The operators will have compact resolvent and therefore the spectrum will consist only of eigenvalues of finite multiplicity.

 Importantly,   the associated energy spaces,  generically denoted by $V(\Omega_\eps)$, will satisfy the condition $ W^{m,2}_0(\Omega_\eps)\subset V(\Omega_{\eps})\subset  W^{m,2}(\Omega_\eps)$ and will depend on the domain and the boundary conditions considered. We will consider different types of boundary conditions according to the choice of the spaces $V(\Omega_\eps)$. If $V(\Omega_\eps)=W^{m,2}_0(\Omega_\eps)$ they will be called Dirichlet boundary conditions, if  $V(\Omega_\eps)=W^{m,2}(\Omega_\eps)$ they will be called Neumann boundary conditions,  and in case $V(\Omega_\eps)=W^{m,2}(\Omega_\eps)\cap W^{k,2}_0(\Omega_\eps)$ for certain $1\leq k\leq m-1$, they will be called ``intermediate boundary condition''. We refer to \cite{Babuska} and references therein for a pioneer discussion on the stability properties under the three types of boundary conditions, including an analysis of the  so-called Babu\v{s}ka-Sapondzhyan paradox.  We also refer to  \cite{Maz-Naz} for a further discussion on the paradox and \cite{Maz-Naz-Pla} for a general reference in this type of problems.  We mention that sharp stability estimates for the eigenvalues of higher order operators subject to Dirichlet and Neumann boundary conditions have been recently proved in \cite{bula2012} where uniform classes of domain perturbations have been considered (see also \cite{bulahigh} for related results); moreover,  in \cite{buoso,buosoplates} further restrictions on the classes of open sets allow obtaining also analyticity results. 

Our goal is twofold.  On one hand, we will provide a rather general condition describing the way the domains converge to the limiting one, which will guarantee the spectral convergence of the operator in $\Omega_\eps$ to the appropriate limiting operator  in $\Omega$. The condition, which we will denote by (C), see Section \ref{stability} below,  is expressed intrinsically and it is posed independently of the boundary conditions that we consider. Needless to say that for a particular family of perturbed domains to check that the condition is satisfied will depend heavily on the boundary condition imposed.   This condition generalizes previous ones formulated for Dirichlet and Neumannn boundary conditions.   

On the other hand, we will focus on the case of higher order operators with ``intermediate boundary'' conditions, paying special attention to the case of the biharmonic operator.  We will obtain almost sharp conditions on the way the boundaries can be perturbed to guarantee the spectral  convergence  with  preservation of the same intermediate boundary conditions for general higher order operators.  Afterwards we will analyze in detail the case of the biharmonic when the boundary of the domain presents an oscillatory behavior. We will see that there is a critical oscillatory behavior such that, when the oscillations are below this threshold  we have  spectral stability, while for oscillations above this value we approach a problem with Dirichlet boundary condition.  For exactly the threshold value, there appears an extra term in the boundary condition for the limiting problem, which maybe interpreted as a ``strange curvature''.   The existence of this critical value is well known in other situations. See for instance the seminal paper \cite{ciora-murat} and also, \cite{Maz-Naz}.   In other context see \cite{casado10, arrbrus2010, komo}.

We describe now the contents of the paper.   In Section \ref{preliminaries} we set up the operators, fix the notation and include a subsection where we describe the basic elements of the technique called ``compact convergence of operators'' which will be used in this paper.  In Section \ref{stability} we state condition (C), see Definition \ref{conc}, and show that this condition implies the  ``compact convergence of operators''  and therefore, the spectral stability of the operators.  In Section \ref{dirichlet} we consider the case of Dirichlet boundary conditions while in Section \ref{neumann} the case of Neumann boundary conditions. The case of intermediate boundary conditions is studied in Section \ref{intermediate}. We prove Lemma \ref{weaklemint} and Corollary \ref{weaklemint2} which provides conditions guaranteeing the spectral stability for the intermediate boundary conditions. 
 These sections cover the first goal of the paper.  

The second goal is achieved in sections \ref{biharmonic-intermediate} and \ref{critical} .  Notice that in Section \ref{biharmonic-intermediate} we analyze the case of a biharmonic operator with intermediate boundary conditions (in this case $V(\Omega)=W^{2,2}(\Omega)\cap W^{1,2}_0(\Omega)$) where the domain is perturbed in an oscillatory way. As a matter of fact if the boundary of the unperturbed domain is given locally around certain point by the function $x_N=g(x_1,\ldots, x_{N-1})$ for $\bar x=(x_1,\ldots, x_{N-1})\in W$ for some nice $(N-1)$-dimensional domain $W$, with $g\equiv 0$ (that is, the boundary is flat) and the boundary of the perturbed domain is given as $x_N=g_\eps(\bar x)$ where $g_\eps(\bar x)=\eps^{\alpha} b(\bar x/\eps)$ for some smooth and periodic function $b$, then $\alpha=3/2$ is a critical value.  If $\alpha>3/2$, the oscillations are not too strong and the limit problem has also the same intermediate boundary conditions. If $\alpha<3/2$, the oscillations are too wild and the limit problem has a Dirichlet boundary condition in $W$. The critical case $\alpha=3/2$ is treated in detail in Section \ref{critical}. We need to treat this case as a homogenization problem and will use the unfolding operator method to show that the limit boundary condition in $W$ contains an extra term. The results of this paper were announced in \cite{ArrLam}.

\section{Preliminaries, notation and some examples}
\label{preliminaries}

%

We introduce in this section the general setting of the paper, the basic notation to follow the contents  and some relevant examples.  Also, we include the definition of ``compact convergence".  

\subsection{Higher order elliptic operators}

We fix some notation and we recall basic facts from standard  spectral theory for elliptic operators. We refer
to Davies~\cite{da} for details and proofs.

Let $N, m \in \mathbb{N}$ and $\Omega$ be an open set in $\mathbb{R}^N$.
We denote by $W^{m,2}(\Omega )$ the Sobolev space  of
real-valued functions in $L^2(\Omega )$, which
have distributional derivatives of order $m$ in $L^2(\Omega ) $,
endowed with the norm
\begin{equation}
\label{norm}
\|u \|_{W^{m,2}(\Omega )}^2=\| u\|_{L^2(\Omega )}^2+ \sum_{|\alpha|= m}\|D^{\alpha }u \|_{L^{2}(\Omega )}^2.
\end{equation}
We denote by  $W^{m,2}_0(\Omega )$
the closure in $W^{m,2}(\Omega )$ of the space of the $C^{\infty}$-functions
with compact support in $\Omega$.

Let $\hat m$ be the number of the
multi-indices $\alpha=(\alpha_1,\dots ,\alpha_N)\in {\mathbb{N}}_0^N$
with length $|\alpha | = \alpha_1+\dots +\alpha_N$ equal to $m$. Here ${\mathbb{N}}_0={\mathbb{N}}\cup\{0\} $.
For all $\alpha ,\beta \in
{\mathbb{N}}_0^N$
such that $|\alpha |=|\beta |=m$, let
$A_{\alpha \beta }$ be bounded measurable real-valued
functions defined on ${\mathbb{R}}^N$ satisfying
$A_{\alpha \beta }=
A_{\beta \alpha }$ and the
condition
\begin{equation}
\label{elp}
\sum_{|\alpha |=|\beta |=m}A_{\alpha \beta }(x)\xi_{\alpha }\xi_{\beta }
\geq 0
\end{equation}
for all $x\in \mathbb{R}^N $, $\xi =(\xi_{\alpha })_{|\alpha |=m}\in\mathbb{R}^{\hat m}$.
For all open sets $\Omega$ in ${\mathbb{R}}^N$ we define  
\begin{equation}\label{quaddef}
Q_{\Omega }(u,v)=\sum_{|\alpha |=|\beta |=m}\int_{\Omega }A_{\alpha \beta }D^{\alpha }uD^{\beta } v\, dx +\int_{\Omega }u v\,dx ,
\end{equation}
for all $u,v\in W^{m,2}(\Omega )$ and
$$Q_{\Omega }(u)=Q_{\Omega }(u,u).$$ Observe that
by condition (\ref{elp})   $Q_{\Omega }$ is in fact a scalar product in $W^{m,2}(\Omega )$.

Let $V(\Omega )$ be a linear subspace of $W^{m,2}(\Omega )$ containing $W^{m,2}_0(\Omega )$. We recall that if
$V(\Omega )$ endowed with the norm $Q_{\Omega }^{1/2}(\cdot )$ is complete then there exists
a uniquely determined non-negative selfadjoint operator $H_{V(\Omega )}$ such that
${\rm Dom }\, H^{1/2}_{V(\Omega )} = V(\Omega )$ and
$$
Q_{\Omega }(u,v)=<H^{1/2}_{V(\Omega )}u, H^{1/2}_{V(\Omega )}u>_{L^2(\Omega )},\ \ \ \forall \ u,v\in V(\Omega ).
$$
In particular, a function $u$ belongs to the domain of $H_{V(\Omega )}$  if and only if  $u\in V(\Omega )$
and there exists $f\in L^2(\Omega )$ such that
\begin{equation}
\label{weak}
Q_{\Omega }(u,v)=<f, v>_{L^2(\Omega )},\ \ \ \forall v\in V(\Omega ).
\end{equation}
Clearly,  $H_{V(\Omega )}u=f$. Equation (\ref{weak}) is the weak formulation of the classical problem
\begin{equation}
\label{classic}
Lu=f,\
\ \ {\rm in}\ \Omega
\end{equation}
where $L$ is the classical operator defined as 
\begin{equation}
\label{classicop}
Lu := (-1)^m\sum_{|\alpha |=|\beta |=m} D^{\alpha }\left(A_{\alpha \beta}D^{\beta }u  \right)+u
\end{equation}
and the unknown  $u$ is subject to suitable homogeneous boundary conditions depending on the choice of $V(\Omega )$ (see the examples below).

We recall that if the embedding $V(\Omega )\subset L^2(\Omega )$ is compact then
the operator $H_{V(\Omega )}$ has compact resolvent. In this case the spectrum is discrete and consists of a sequence of eigenvalues $\lambda_n[V(\Omega )]$
of finite multiplicity which can be represented by means of the Min-Max Principle:
$$
\lambda_n[V(\Omega )]= \inf_{\substack{E\subset V(\Omega )\\ {\rm dim }\, E=n  }}\sup_{\substack{u\in E\\ u\ne 0}}\frac{Q_{\Omega }(u)}{\| u \|^2_{L^2(\Omega )}}.
$$
Correspondingly, there exists an orthonormal basis in $L^2(\Omega )$ of eigenfunctions $\varphi_n[V(\Omega )]$ associated with the eigenvalues $\lambda_n[V(\Omega )]$.

Note that, since the coefficients $A_{\alpha \beta }$ are fixed and bounded then
$$
Q_{\Omega }^{1/2}(u)\le C \| u\|_{W^{m,2}(\Omega )},
$$
for all $u\in W^{m,2}(\Omega )$ where $C$ is a positive constant independent of $u$ and $\Omega$. Thus, since we assume that the space $(V(\Omega ), Q_{\Omega }^{1/2}(\cdot ))$ is complete, we have that
$$
  c\| u\|_{W^{m,2}(\Omega )} \le Q_{\Omega }^{1/2}(u),
$$
for all $u\in V(\Omega )$, where $c$ is a positive constant independent of $u$. In other words, the two norms $Q_{\Omega }^{1/2}(\cdot )$
and $\| \cdot \|_{W^{m,2}(\Omega )}$ are equivalent in $V(\Omega )$. Note that in general the constant $c$ may depend on $\Omega$. However, stronger assumptions on the coefficients allow us to get $c$ independent of $\Omega$. For example, if
the coefficients $A_{\alpha  \beta }$ satisfy the uniform ellipticity condition
\begin{equation}
\label{unielp}
\sum_{|\alpha |=|\beta |=m}A_{\alpha \beta }(x)\xi_{\alpha }\xi_{\beta }
\geq \theta \sum_{|\alpha  |=m}|\xi_{\alpha }|^2
\end{equation}
for all $x\in \mathbb{R}^N $, $\xi =(\xi_{\alpha })_{|\alpha |=m}\in\mathbb{R}^{\hat m}$, then it is straightforward that  $c$ can be chosen
$c=\min\{\sqrt{\theta},1\}$ which is independent of $\Omega$. Condition (\ref{unielp}) will not be used in Section 2 which is devoted to a general stability theorem. However, we shall use it in the following sections devoted to applications.

\subsection{${\mathcal{E}}$-compact convergence}

Let $\Omega$ be a fixed open set and $V(\Omega )$ its corresponding space as in the previous section. For all sufficiently small $\epsilon >0$ we consider perturbations $\Omega_{\epsilon }$ of $\Omega$ and we denote by $V(\Omega_\epsilon)$ the corresponding spaces of functions defined on $\Omega_{\epsilon }$. 
We assume that the coefficients $A_{\alpha  \beta }$ are fixed functions defined in the whole of ${\mathbb{R}}^N$ and that the operators $H_{V(\Omega )}$ and $H_{V(\Omega_{\epsilon } )}$ are well-defined and have compact resolvents.

We denote by ${\mathcal E}$ the extension-by-zero operator, which means that given a real-valued function $u$ defined on some set in ${\mathbb{R}}^N$, ${\mathcal E}u$ is the function extended by zero outside the given set.  Clearly, for each $\epsilon >0$,  ${\mathcal E}$ can be thought as an operator  acting from $L^2(\Omega )$ to $L^2(\Omega _{\epsilon })$, consisting in extending the function by zero to all of $\mathbb{R}^N$ and then restricting it to $\Omega_\eps$.  As a matter of fact, this operator ${\mathcal E}$ will be the key to compare functions and operators defined in $\Omega$ and $\Omega_\eps$.  The following concepts and definitions go back  to the works of F. Stummel, see \cite{Stu} and G. Vainniko see \cite{Vai77a,Vai79} among others.  We also refer to \cite{CarPis,ArrCarLoz}.    Here we denote by ${\mathcal{L}}(X)$ the space of bounded linear operators acting from a normed  space $X$ to itself.

\begin{defn}\label{E-comp-conv}
i) We say that $v_\eps\in L^2(\Omega_\eps)$ $\mathcal{E}$-converges to $v\in L^2(\Omega)$ if $\|v_\eps-\mathcal{E}v\|_{L^2(\Omega_\eps)}$ $\to 0$ as $\eps\to 0$. We write this as 
$v_\eps\Econ v$.
\par\noindent ii) The family of bounded linear operators $B_\eps\in \mathcal{L}(L^2(\Omega_\eps))$ $\mathcal{EE}$- converges to $B\in \mathcal{L}(L^2(\Omega))$ if $B_\eps v_\eps\Econ Bv$ whenever $v_\eps\Econ v$. We write this as 
$B_\eps\EEcon B$.

\par\noindent iii) The family of bounded linear and compact  operators $B_\eps\in \mathcal{L}(L^2(\Omega_\eps))$ $\mathcal{E}$-compact converges to $B\in \mathcal{L}(L^2(\Omega))$ if $B_\eps\EEcon B$ and for any family of functions $v_\eps\in L^2(\Omega_\eps)$ with $\|v_\eps\|_{L^2(\Omega_\eps)}\leq 1$ there exists a subsequence, denoted by $v_\eps$ again, and a function $w\in L^2(\Omega)$ such that $B_\eps v_\eps \Econ w$. We write
$B_\eps\Ccon B$.
\end{defn} 

\par\medskip 

\par\medskip
There is a strong relation between the  $\mathcal{E}$-compact convergence of a family of operators and their spectral convergence.  By  this, we mean the  convergence of eigenvalues and the associated spectral projections, see \cite[Section 2.1]{arca}. Since in this particular work  we are mainly dealing with $B$ and $B_\eps$ which are the inverses of the operators $H_{V(\Omega)}$ and $H_{V(\Omega_\eps)}$ defined above, we will define the spectral convergence just for this special type of operators. Hence, if we denote by $\{(\lambda_n^\eps,\phi_n^\eps)\}_{n=1}^\infty$ the eigenvalues and eigenfunctions of $H_{V(\Omega_\eps)}$ and by $\{(\lambda_n,\phi_n)\}_{n=1}^\infty$ the eigenvalues and eigenfunctions of $H_{V(\Omega)}$, (where we understand that the eigenfunctions are extended by zero outside $\Omega_\eps$ and $\Omega$ and they are normalized in $L^2(\mathbb{R}^N)$, we will say that the spectra behaves continuously at $\eps=0$,
if for fixed $n\in {\mathbb N}$ we have that $\la_n^\eps\to \la_n$ as $\eps\to 0$ and the
spectral projections converge in $L^2(\mathbb{R}^N)$, that is, if $a\not\in
\{\la_n \}_{n=1} ^\infty$, and $\la_n <a<\la_{n+1} $,
 then if we define the projections
$P_a^\eps:L^2(\mathbb{R}^N)\to L^2(\mathbb{R}^N)$, $P_a^\eps(\psi)=\sum_{i=1}^n(\phi_i
^\eps,\psi)_{L^2(\mathbb{R}^N)}\phi_i^\eps$ then
$$\sup\{\|P_a^\eps(\psi)-P_a^0(\psi)\|_{L^2(\mathbb{R}^N)}:\quad \psi\in
L^2(\mathbb{R}^N),\, \|\psi\|_{L^2(\mathbb{R}^N)}=1\}\to
0,\quad \hbox{as }\eps\to 0.$$

The convergence of the spectral projections is equivalent to
the following: for each sequence $\eps_k\to 0$ there exists a
subsequence, that we denote again by $\eps_k$, and a complete
system of orthonormal eigenfunctions of the limiting problem
$\{\phi_n\}_{n=1}^\infty$ such that
$\|\phi_n^{\eps_k}-\phi_n\|_{L^2(\mathbb{R}^N)}\to 0$ as $k\to \infty$.


As a matter of fact, we can show the following (see \cite[Thm.~4.1]{ArrCarLoz}, \cite{CarPis}).  

\begin{prop}\label{spectral-convergence} Assume the operator $\mathcal{E}$ satisfies the condition $\|\mathcal{E} u\|_{L^2(\Omega_\eps)}\to \|u\|_{L^2(\Omega)}$ for each $u\in L^2(\Omega)$.  If $H_{V(\Omega_\eps)}^{-1}\in \mathcal{L}(L^2(\Omega_\eps))$ are compact and $H_{V(\Omega_\eps)}^{-1} \Ccon H_{V(\Omega)}^{-1} $, then  we have the spectral convergence of $H_{V(\Omega_\eps)}$ to $H_{V(\Omega)}$.   
\end{prop}

\subsection{Examples}
We consider in this section some relevant examples of higher order operators. 

\subsubsection{Polyhamonic operators}

An important class of higher order operators is given by the polyharmonic operators which we brefly discuss here as a prototypical example,   see for instance \cite{gaz}. 

For  $m\in {\mathbb{N}}$, we set   $A_{\alpha\beta}= \delta_{\alpha \beta }m!/ \alpha ! $ for all $\alpha, \beta \in {\mathbb{N}}^N$ with $|\alpha |=|\beta |=m$, where $\delta _{\alpha \beta }=1$ if $\alpha =\beta $ and $\delta_{\alpha\beta}=0$ otherwise.  With this choice, condition (\ref{unielp})
is satisfied. 

%

Let $k\in {\mathbb{N}}_0$, $0\le k \le m $ and $V(\Omega )=W^{m,2}(\Omega )\cap W^{k,2}_0(\Omega )$, endowed with the norm (\ref{norm}) of $W^{m,2}(\Omega)$. If $k=m$ and $\Omega $ has finite Lebesgue measure then $V(\Omega )$ is a closed subspace of $W^{m,2}(\Omega )$ and the embedding $V(\Omega )\hookrightarrow L^2(\Omega )$ is compact.
If $0<k<m$ then, under very weak regularity assumptions on $\Omega$ (say, $\Omega $ has a quasi-resolved boundary\footnote{this includes the case of open sets satisfying the classical cone condition, as well as the case of open sets  of class $C^0$, see Definition~\ref{atlas} below.} in the sense of Burenkov\cite[\S 4.3]{bu}), $V(\Omega )$ is a closed subspace of $W^{m,2}(\Omega)$ and $V(\Omega )\hookrightarrow W^{k,2}_0(\Omega )$; if in addition $\Omega $ has finite Lebesgue measure then the embedding $W^{k,2}_0(\Omega )\hookrightarrow L^2(\Omega )$ is compact, hence the embedding $V(\Omega )\hookrightarrow L^2(\Omega )$ is also compact.
If $k=0$  and the open set $\Omega $ is bounded and  of class $C^0$, see  Definition~\ref{atlas} below, then the embedding $V(\Omega )\hookrightarrow L^2(\Omega )$
is compact (see Burenkov~\cite[Thm.~8,~p.169]{bu}).

Note that if $k=m$ then $V(\Omega )= W^{m,2}_0(\Omega )$ and
integrating by parts one can  realize that
$$\sum_{|\alpha |=|\beta |=m}\int_{\Omega }A_{\alpha \beta }D^{\alpha }uD^{\beta } v\, dx 
=\left\{\begin{array}{ll}\int_{\Omega }\Delta^{\frac{m}{2}}u\Delta^{\frac{m}{2}}v dx,\ \ &{\rm if }\ m\ {\rm is\ even\, ,}\vspace{2mm}\\
\int_{\Omega }\nabla \Delta^{\frac{m-1}{2}}u\nabla \Delta^{\frac{m-1}{2}}v dx,\ \ &{\rm if }\ m\ {\rm is\ odd\, , }
\end{array}\right.
$$
for all $u,v \in W^{m,2}_0(\Omega )$. In this case we obtain in (\ref{classicop}) the  operator,   $Lu=(-\Delta)^m u+u$  subject to
the Dirichlet boundary conditions  $u=\frac{\partial u}{\partial \nu}=\dots \frac{\partial^{m-1} u}{\partial \nu^{m-1}}=0$ on $\partial \Omega$.  Here and in the sequel $\nu $ denotes the unit outer normal to $\partial \Omega$. The operator $(-\Delta)^m u$ is the classical polyharmonic operator of order $2m$.

In the general case $k\le m$, the classical problem  reads
$$\left\{\begin{array}{ll}
(-\Delta )^mu +u=f,\ \ & {\rm in }\ \Omega,\\
\frac{\partial^j u}{\partial \nu^j}=0, \ \forall\ j=0,\dots , k-1,\ \ & {\rm on }\ \partial \Omega , \\
{\mathcal{B}}_ju=0, \ \forall\ j=1,\dots , m-k, \ \ & {\rm on }\ \partial \Omega ,
\end{array}\right. $$
where ${\mathcal{B}}_j$ are uniquely defined `complementing' boundary operators. See Nec\v{a}s~\cite{nec} for details.  \\

\subsubsection{Biharmonic operator with Dirichlet, Neumann and Intermediate boundary conditions}

Let us consider the case  $m=2$ in the previous example. The quadratic form (\ref{quaddef}) will read 
\begin{equation}
\int_{\Omega }D^2u:D^2v dx +\int_{\Omega}uvdx, 
\end{equation}
for all $u, v$ in an appropriate energy space $V(\Omega)$. Here and in the sequel  $D^2u$ denotes the Hessian
matrix of $u$ and $D^2u: D^2v=\sum_{i,j=1}^N\frac{\partial^2u}{\partial x_i\partial x_j}\frac{\partial^2v}{\partial x_i\partial x_j}$ is the Frobenius product of the two matrices.

As above, if $V(\Omega)=W^{2,2}_0(\Omega)$,   the classical operator (\ref{classic}) is given   by the biharmonic operator and we obtain the classical Dirichlet  problem 
$$\left\{\begin{array}{ll}
\Delta^ 2u +u=f,\ \ & {\rm in }\ \Omega,\\
u=0,\ \ & {\rm on }\ \partial \Omega , \\
\frac{\partial u}{\partial \nu}=0,\ \ & {\rm on }\ \partial \Omega . \\
\end{array}\right. $$

It is well-known that if $N=2$ the Dirichlet problem for the biharmonic operator is related  for example to  the study of the bending of clamped plates. 

If $V(\Omega)= W^{2,2}(\Omega)$, by using the `Biharmonic Green Formula' (\ref{Green}) we obtain the classical Neumann problem

$$\left\{\begin{array}{ll}
\Delta^2u +u=f,\ \ & {\rm in }\ \Omega,\\
\frac{\partial ^2u}{\partial \nu ^2}=0,\ \ & {\rm on }\ \partial \Omega , \\
{\rm div}_{\partial  \Omega }((D^2 u)\cdot \nu)_{\partial \Omega}  +\frac{\partial \Delta u}{\partial \nu}=0,\ \ & {\rm on }\ \partial \Omega ,\\
\end{array}\right. $$
involving the well-known tangential divergence operator, see Section~\ref{microscopic} for basic definitions.   
We recall that if $N=2$ the Neumann problem for the biharmonic operator arises for example in the study of the bending of free plates. See also Chasman~\cite{chas}.

Finally,  If $V(\Omega)= W^{2,2}(\Omega)\cap W^{1,2}_0(\Omega)$, proceeding as above we obtain the classical intermediate  problem

\begin{equation}\label{simplysupp}\left\{\begin{array}{ll}
\Delta^2u +u=f,\ \ & {\rm in }\ \Omega,\\
u=0,\ \ & {\rm on }\ \partial \Omega , \\
\frac{\partial ^2u}{\partial \nu ^2}=0,\ \ & {\rm on }\ \partial \Omega . \\
\end{array}\right. 
\end{equation}

We recall that if $N=2$ the intermediate  problem for the biharmonic operator arises for example in the study of the bending of hinged plates (sometimes called simply-supported).

We note that since $u=0$ on $\partial \Omega $ then the second boundary condition in (\ref{simplysupp}) can be written as 
$$
\Delta u-K\frac{\partial u}{\partial \nu }=0,
$$
where $K$ is the mean curvature of the boundary, i.e., the sum of the principal curvatures. 

See Gazzola, Grunau and Sweers~\cite{gaz} for further details.

\section{A general stability theorem }
\label{stability}

%

The following condition on  $\Omega_{\epsilon }$ and $V(\Omega_{\epsilon })$ will guarantee that
$H_{V(\Omega_{\epsilon })}^{-1}$  is ${\mathcal{E}}$-compact convergent  to $H_{V(\Omega )}^{-1}$  in the sense of Definition~\ref{E-comp-conv}.\\

\begin{defn}{\bf (Condition C)}\label{conc}
 Given open sets $\Omega_{\epsilon}$, $\epsilon>0$ and $\Omega $ in ${\mathbb{R}}^N$ and corresponding elliptic operators
$ H_{V(\Omega_{\epsilon })}$, $H_{V(\Omega)}$ defined on $\Omega_{\epsilon}, \Omega $ respectively, we say that condition {\bf (C)} is satisfied if
for each $\epsilon >0$ there exists an open set $K_{\epsilon }\subset \Omega \cap \Omega_{\epsilon }$ such that\footnote{this condition guarantees that $\| {\mathcal{E}}u \|_{L^2(\Omega_{\epsilon })}\to \| u \|_{L^2(\Omega )}$ which is the basic hypothesis in the theory of $E$-convergence and it is an assumption of Proposition~\ref{spectral-convergence}. }
\begin{equation}\label{limmeas}\lim_{\epsilon \to 0}| \Omega \setminus K_{\epsilon }  |= 0, \end{equation}
 and the following conditions are satisfied:\\

\noindent {\bf (C1)} If $v_{\epsilon }\in V(\Omega_{\epsilon })$ and $\sup_{\epsilon >0}Q_{\Omega_{\epsilon }}(v_{\epsilon })<\infty $ then
\begin{equation}\label{c1lim} \lim_{\epsilon \to 0} \| v_{\epsilon }\|_{L^2(\Omega_{\epsilon }\setminus K_{\epsilon })} =0;\end{equation}
\noindent {\bf (C2)}
For each $\epsilon >0$ there exists an operator
$T_{\epsilon }:V(\Omega )\longrightarrow V(\Omega_{\epsilon })$
such that for all fixed $\varphi\in V(\Omega)$
\begin{itemize}
\item[{\rm (i)}]\  $\displaystyle\lim_{\epsilon\to 0}Q_{K_{\epsilon }}(T_{\epsilon } \varphi -\varphi)=0,$
\item[{\rm (ii)}]\  $\displaystyle\lim_{\epsilon\to 0}Q_{\Omega_{\epsilon }\setminus K_{\epsilon }}( T_{\epsilon} \varphi )=0,$
\item[{\rm (iii)}]\  $\displaystyle\sup_{\epsilon >0}\| T_{\epsilon }\varphi \|_{L^2(\Omega _{\epsilon })} <\infty .$
\end{itemize}

\noindent   {\bf  (C3)}
For each $\epsilon >0$ there exists an operator
$E_{\epsilon }$ from $V(\Omega_{\epsilon})$ to $W^{m,2}(\Omega )$ such that the set $E_{\eps}(V(\Omega_{\epsilon}))$ is compactly embedded in $L^2(\Omega )$ and such that  
\begin{itemize}
\item[{\rm (i)}] If $v_{\epsilon}\in V(\Omega_{\epsilon})$ is such that $\sup_{\epsilon >0}Q_{\Omega_{\epsilon} }(v_{\epsilon })<\infty $  then  $\displaystyle\lim_{\epsilon\to 0}Q_{K_{\epsilon }}(E_{\epsilon } v_{\epsilon } -v_{\epsilon })=0$;
\item[{\rm (ii)}]\  $\displaystyle\sup_{\epsilon >0}\sup_{\substack{v \in V(\Omega_{\epsilon })\\ v \ne 0}}\frac{ \| E_{\epsilon }v \|_{W^{m,2}(\Omega )} }{Q_{\Omega_{\epsilon} }^{1/2}(v )}<\infty $;
\item[{\rm (iii)}]\  If $v_{\epsilon}\in V(\Omega_{\epsilon})$ is such that $\sup_{\epsilon >0}Q_{\Omega_{\epsilon} }(v_{\epsilon })<\infty $ and
there exists $v \in L^2(\Omega )$ such that, possibly passing to a subsequence, we have 
$\|E_\eps v_\eps-v\|_{L^2(\Omega)}\eto 0$, then $v \in V(\Omega )$.

\end{itemize}
\end{defn}

\begin{example} Consider the simpler case   $\Omega\subset \Omega_{\epsilon }$  and $V(\Omega)=W^{m,2}(\Omega)$, $V(\Omega_\eps)=W^{m,2}(\Omega_\eps)$  for all $\epsilon >0$, and assume  (\ref{unielp}) is satisfied. We set $K_{\epsilon}=\Omega $ for all $\epsilon >0$. Thus (\ref{limmeas}) is trivially satisfied.  Assume that $\Omega $ is sufficiently regular to guarantee the existence  of a bounded extension operator   from $W^{m,2}(\Omega )$ to $W^{m,2}({\mathbb{R}}^N )$. Then condition {\bf (C2)} is satisfied: indeed, the extension operator may serve as operator  $T_{\epsilon }$. As far as the operator $E_{\epsilon }$ is concerned,  one  can use the restriction operator: in this  way,  the compactness of the embedding 
$W^{m,2}(\Omega )\hookrightarrow L^2(\Omega)$  allow to conclude that condition {\bf (C3)} is satisfied . 
Thus,  in order to verify  the validity of condition {\bf (C)} it just suffices to  check  the validity of {\bf (C1)}.

%
%
\end{example}

We now prove the following general statement.

\begin{thm}\label{general} Let ${\mathcal{E}}$ be the extension-by-zero operator.
If condition (C) is satisfied then $H^{-1}_{\Omega_{\epsilon}}\Ccon  H^{-1}_{\Omega}$.
\end{thm}

\par\noindent {\bf Proof.} First we prove that $H^{-1}_{\Omega_{\epsilon}}\EEcon H^{-1}_{\Omega}$.
Let $f_{\epsilon }\in L^2(\Omega_{\epsilon })$ be a sequence ${\mathcal {E}}$-convergent to $f\in L^2(\Omega )$, {\it i.e.,}
\begin{equation}
\label{e05}\lim_{\epsilon \to 0}
\|f_{\epsilon }-{\mathcal{E}}f\|_{L^2(\Omega_{\epsilon })}=0
.\end{equation}
Let $v_{\epsilon }\in V(\Omega _{\epsilon })$ be such that
\begin{equation}
\label{e1}
H_{\Omega_{\epsilon}}v_{\epsilon }=f_{\epsilon }  ,
\end{equation}
 and let $u\in V(\Omega )$ be such that $H_{\Omega }u=f$. We have to prove that
$v_{\epsilon }$ is ${\mathcal {E}}$-convergent to $u$, {\it i.e.,}
\begin{equation}
\label{e1.5}
\lim_{\epsilon \to 0}\| v_{\epsilon }-{\mathcal {E}} u\|_{L^2(\Omega _{\epsilon })}\to 0 .
\end{equation}
 We will prove this statement by showing that for any sequence $\eps_k\to 0$ there is a subsequence $\eps_k'\to 0$ which satisfies \eqref{e1.5}. Moreover, in order to avoid a complicated notation with too many indices and subindices, we will keep denoting the sequences and subsequences by $\eps$.  

By (\ref{e1}) and the H\"{o}lder inequality   it follows that
\begin{equation}
\label{e2}
Q_{\Omega_{\epsilon }}(v_{\epsilon })=(f_{\epsilon },v_{\epsilon })_{L^2(\Omega_{\epsilon })}\le \| f_{\epsilon }\|_{L^2(\Omega_{\epsilon })}Q^{1/2}_{\Omega_{\epsilon }}(v_{\epsilon }).
\end{equation}
Since $f_{\epsilon }$ is ${\mathcal{E}}$-convergent to $f$, there exists $M>0$ such that $\| f_{\epsilon }\|_{L^2(\Omega_{\epsilon })}\le M$ for all $\epsilon >0$. Thus by (\ref{e2}) it follows that
\begin{equation}
\label{e3}
\sup_{\epsilon >0}Q_{\Omega_{\epsilon }}(v_{\epsilon })<\infty .
\end{equation}
By condition (C3) (ii), it follows that $E_{\epsilon }v_{\epsilon }$ is a bounded sequence in $W^{m,2}(\Omega )$. Accordingly, by the compactness of the embedding 
$E_{\eps}(V(\Omega_{\eps}))\subset L^2(\Omega )$ and the reflexivity of the space $W^{m,2}(\Omega )$ there exists
$\tilde u\in W^{m,2}(\Omega )$ such that,
possibly considering a subsequence, $E_{\epsilon }v_{\epsilon }$ converges to $\tilde u$ strongly in $L^2(\Omega  )$ and weakly in $W^{m,2}(\Omega )$ as $\epsilon \to 0$. By (\ref{e3}) and condition (C3) (iii) it follows that $\tilde u \in V(\Omega )$. We now prove that $\tilde u=u$.
Let $\varphi \in V(\Omega )$ be fixed. Since the operator $T_{\epsilon }$ takes values in $V(\Omega _{\epsilon })$ we can use $T_{\epsilon}\varphi $ as a test
function in the weak formulation of the problem in $\Omega_{\epsilon }$ and obtain
\begin{equation}
\label{e4}
Q_{\Omega_{\epsilon }}(v_{\epsilon }, T_{\epsilon }\varphi )=(f_{\epsilon }, T_{\epsilon }\varphi ) ,
\end{equation}
It is easily seen that
\begin{eqnarray}
\lefteqn{ Q_{\Omega_{\epsilon }}(v_{\epsilon }, T_{\epsilon }\varphi ) =
Q_{K_{\epsilon }}(v_{\epsilon }, T_{\epsilon }\varphi )+Q_{\Omega_{\epsilon }\setminus K_{\epsilon }}(v_{\epsilon }, T_{\epsilon }\varphi )
}\nonumber \\ & &\qquad
   =Q_{K_{\epsilon }}(v_{\epsilon }, \varphi )+ Q_{K_{\epsilon }}(v_{\epsilon }, T_{\epsilon }\varphi -\varphi  ) +Q_{\Omega_{\epsilon }\setminus K_{\epsilon }}(v_{\epsilon }, T_{\epsilon }\varphi )   \nonumber \\
& & \qquad  =  Q_{K_{\epsilon }}(E_{\epsilon }v_{\epsilon }, \varphi )+Q_{K_{\epsilon }}(v_{\epsilon }-E_{\epsilon }v_{\epsilon }, \varphi )
+ Q_{K_{\epsilon }}(v_{\epsilon }, T_{\epsilon }\varphi -\varphi ) +Q_{\Omega_{\epsilon }\setminus K_{\epsilon }}(v_{\epsilon }, T_{\epsilon }\varphi )   \nonumber \\
& & \qquad  =
Q_{\Omega   }(E_{\epsilon }v_{\epsilon }, \varphi )- Q_{\Omega \setminus K_{\epsilon  }   }(E_{\epsilon }v_{\epsilon }, \varphi )
+ Q_{K_{\epsilon }}(v_{\epsilon }-E_{\epsilon }v_{\epsilon }, \varphi )
+ Q_{K_{\epsilon }}(v_{\epsilon }, T_{\epsilon }\varphi -\varphi )
 \nonumber \\
& &\qquad \ \ \ +Q_{\Omega_{\epsilon }\setminus K_{\epsilon }}(v_{\epsilon }, T_{\epsilon }\varphi ).
\label{e5}
\end{eqnarray}

By the boundedness of the coefficients $A_{\alpha \beta }$, the space $(W^{m,2}(\Omega ), \|\cdot \|_{W^{m,2}(\Omega )})$  is continuously embedded
in $(W^{m,2}(\Omega ), Q^{1/2}_{\Omega }(\cdot ))$. It follows that the sequence $E_{\epsilon }v_{\epsilon }$ is weakly convergent to
$\tilde u $ in $(W^{m,2}(\Omega ), Q^{1/2}_{\Omega }(\cdot ))$, hence
\begin{equation}
\label{e6}
\lim_{\epsilon \to 0}Q_{\Omega }(E_{\epsilon }v_{\epsilon }, \varphi  )=Q_{\Omega }(\tilde u , \varphi ).
\end{equation}

Since $Q_{\Omega \setminus K_{\epsilon }}(E_{\epsilon }v_{\epsilon })$ is a bounded sequence, $\varphi $ is fixed and $\lim_{\epsilon \to 0} |\Omega \setminus K_{\epsilon }|=0$ it is easily seen that
\begin{equation}
\label{e9}
\lim_{\epsilon \to 0}Q_{\Omega \setminus K_{\epsilon }}(E_{\epsilon }v_{\epsilon }, \varphi )=0.
\end{equation}

Since $Q_{K_{\epsilon }}(\varphi )$ is a bounded sequence, by condition (C3) (i) it follows that
\begin{equation}
\label{e65}
\lim_{\epsilon \to 0}Q_{ K_{\epsilon }}(v_{\epsilon } -E_{\epsilon }v_{\epsilon }, \varphi )=0.
\end{equation}

By (\ref{e3}) $Q_{K_{\epsilon }}(v_{\epsilon })$ is a bounded sequence hence by (C2) (i) it follows that
\begin{equation}
\label{e8}
\lim_{\epsilon \to 0} Q_{K_{\epsilon }}(v_{\epsilon }, T_{\epsilon }\varphi -\varphi )=0. \nonumber
\end{equation}

Similarly, by (C2) (ii)
\begin{equation}
\label{e10}
\lim_{\epsilon \to 0}  Q_{\Omega_{\epsilon } \setminus K_{\epsilon }}(v_{\epsilon }, T_{\epsilon }\varphi )=0.
\end{equation}

Thus, by (\ref{e5})-(\ref{e10}) it follows that
\begin{equation}
\label{e11}
 \lim_{\epsilon \to 0}Q_{\Omega_{\epsilon }}(v_{\epsilon }, T_{\epsilon }\varphi )=Q_{\Omega }(\tilde u , \varphi )
\end{equation}

Moreover,
\begin{equation}
\label{e12}
(f_{\epsilon }, T_{\epsilon }\varphi )_{L^2(\Omega_{\epsilon })}=(f_{\epsilon }-{\mathcal{E}}f, T_{\epsilon }\varphi )_{L^2(\Omega_{\epsilon })}
+  ({\mathcal{E}}f, T_{\epsilon }\varphi )_{L^2(\Omega_{\epsilon })}.
\end{equation}

By (\ref{e05}) and condition (C2) (iii) we have
\begin{equation}
\label{e13}
\lim_{\epsilon\to 0} (f_{\epsilon }-{\mathcal{E}}f, T_{\epsilon }\varphi )_{L^2(\Omega_{\epsilon })}=0.
\end{equation}

Furthermore,
\begin{eqnarray}
\label{e14}
\lefteqn{
|({\mathcal{E}}f, T_{\epsilon }\varphi )_{L^2(\Omega_{\epsilon })}-(f, \varphi )_{L^2(\Omega )}|
\le
|(f, T_{\epsilon }\varphi -\varphi )_{L^2(K_{\epsilon })}| }\nonumber  \\
& & \qquad\qquad\qquad
 +|(f, T_{\epsilon }\varphi )_{L^2((\Omega \cap \Omega_{\epsilon })\setminus K_{\epsilon})}|
 +   |(f, \varphi )_{L^2(\Omega \setminus K_{\epsilon})}| .
\end{eqnarray}

By condition (C2) (i) the first summond  in the right-hand side of (\ref{e14}) vanishes as $\epsilon\to 0$. Moreover,
since $|(\Omega \cap \Omega_{\epsilon })\setminus K_{\epsilon}|\to 0 $ as $\epsilon \to 0$,  by condition (C2) (iii) the second and the third summonds in the right-hand side of (\ref{e14}) vanish as $\epsilon\to 0$. Thus
\begin{equation}
\label{e15}\lim_{\epsilon\to 0 }
({\mathcal{E}}f, T_{\epsilon }\varphi )_{L^2(\Omega_{\epsilon })}
= (f, \varphi )_{L^2(\Omega )} .
\end{equation}

By (\ref{e4}), (\ref{e11}), (\ref{e12}), (\ref{e15}) it follows that
\begin{equation}
\label{e16}
Q_{\Omega }(\tilde u , \varphi)=
(f, \varphi )_{L^2(\Omega )},\quad \forall \varphi\in V(\Omega).
\end{equation}
Since, $\tilde u\in V(\Omega)$, we have $\tilde u =u$.

Observe now that
\begin{eqnarray}
\| v_{\epsilon}- {\mathcal{E}}u \|^2_{L^2(\Omega_{\epsilon })} & =&  \| v_{\epsilon} \|^2_{L^2(\Omega_{\epsilon }\setminus \Omega )}+\| v_{\epsilon}- u \|^2_{L^2(\Omega_{\epsilon }\cap  \Omega )}\nonumber \\
& =&  \| v_{\epsilon} \|^2_{L^2(\Omega_{\epsilon }\setminus \Omega )}+\| v_{\epsilon}- u \|^2_{L^2( (\Omega_{\epsilon }\cap  \Omega )\setminus K_{\epsilon }) } +\| v_{\epsilon}- u \|^2_{L^2( K_{\epsilon }) }.
\label{e17}
\end{eqnarray}
Note that
\begin{eqnarray}
\label{e18}
\| v_{\epsilon}- u \|_{L^2(K_{\epsilon } )}^2\le \| E_{\epsilon }v_{\epsilon}-  u \|^2_{L^2(K_{\epsilon })}+  \|  E_{\epsilon }v_{\epsilon }-v_{\epsilon } \|^2_{L^2(K_{\epsilon })}\nonumber  \\
\le \| E_{\epsilon }v_{\epsilon}-  u \|^2_{L^2(\Omega )}+  \|  E_{\epsilon }v_{\epsilon }-v_{\epsilon } \|^2_{L^2(K_{\epsilon })}
\end{eqnarray}
Moreover,
\begin{equation}
\label{e195}
\| v_{\epsilon} \|^2_{L^2(\Omega_{\epsilon }\setminus \Omega )}\le \| v_{\epsilon} \|^2_{L^2(\Omega_{\epsilon }\setminus K_{\epsilon } )},
\end{equation}
and
\begin{equation}
\label{e196}
\| v_{\epsilon}- u \|^2_{L^2( (\Omega_{\epsilon }\cap  \Omega )\setminus K_{\epsilon }) }\le \| v_{\epsilon}\|^2_{L^2( \Omega_{\epsilon } \setminus K_{\epsilon }) }+\| u\|^2_{L^2(   \Omega \setminus K_{\epsilon }) }.
\end{equation}
Thus, by (C1), (C3) (i) and (\ref{e17})-(\ref{e196}) it follows that
(\ref{e1.5}) holds.  Thus,  $H^{-1}_{\Omega_{\epsilon}}$ is ${\mathcal{E}}$-convergent to  $ H^{-1}_{\Omega}$ as $\epsilon \to 0$.

Exactly the same argument can be used to prove
that if  $\hat f_{\epsilon }\in L^2(\Omega _{\epsilon })$ is such that
$$
\sup_{\epsilon >0}\| \hat f_{\epsilon }\|_{L^2(\Omega_{\epsilon })}<\infty
$$
and $\hat v_{\epsilon }\in V(\Omega_{\epsilon })$ is such that $H_{\Omega_{\epsilon }}\hat v_{\epsilon }=\hat f_{\epsilon } $ then
there exist $\hat u\in L^2(\Omega )$ such that, possibly considering a subsequence,
\begin{equation}
\label{e21}
\lim_{\epsilon \to 0} \| \hat v_{\epsilon}- {\mathcal{E}}\hat u \|^2_{L^2(\Omega_{\epsilon })}=0.
\end{equation}
 This implies that $H^{-1}_{\Omega_{\epsilon}}$  is ${\mathcal{E}}$-compact  convergent  to  $ H^{-1}_{\Omega}$ as $\epsilon \to 0$. \hfill $\Box$

\section{Dirichlet boundary conditions and  Mosco
 convergence }
\label{dirichlet}

In this section we consider the operator (\ref{classicop}) on a bounded open set $\Omega $ in ${\mathbb{R}}^N$, subject to Dirichlet boundary conditions 
\begin{equation}
\label{dircon}
u=\frac{\partial u}{\partial \nu}=\dots \frac{\partial^{m-1} u}{\partial \nu^{m-1}}=0,\ \ {\rm on }\ \partial \Omega. 
\end{equation}
This has to be understood in the general frame discussed in Section 2 as follows.

Imposing Dirichlet boundary conditions to the operator $L$ on $\Omega $  means that the domain $V(\Omega )$ of the corresponding quadratic form $Q_{\Omega }$ is given by
$$
V(\Omega )=W^{m,2}_0(\Omega ).
$$
This will be understood throughout this section.
Here we assume that the coefficients $A_{\alpha  \beta }$ satisfy the uniform ellipticity condition
(\ref{unielp}).
We recall that if $\Omega $ is bounded then the Sobolev space $W^{m,2}_0(\Omega )$ is compactly embedded in $L^2(\Omega )$. Thus, as it is explained in Section 2, the operator $H_{W^{m,2}_0(\Omega )}$ is well-defined and has compact resolvent.

The spectral stability of higher order operators subject to Dirichlet boundary conditions on variable domains
was discussed in Babuska and Vyborny~\cite{babvyb} where sufficient conditions ensuring stability were given. Those conditions
are nowadays understood in the frame of the notion of Mosco convergence which we now recall.

\begin{defn}\label{Mosco} Let $D$ be a bounded open set in ${\mathbb{R}}^N$. Let $\Omega_{\epsilon}$, $\epsilon >0$ be a family of open sets contained in $D$. Let $m\in {\mathbb{N}}$ and $\Omega $ be an open set in $D$. We say that the spaces $W^{m,2}_0(\Omega_{\epsilon })$ converge in the sense of Mosco to the space $W^{m,2}_0(\Omega )$ as $\epsilon\to 0$ if the following two conditions are satisfied:
\begin{itemize}
\item[{\rm (M1)}] For any $\varphi\in W^{m,2}_0(\Omega )$ and $\epsilon >0$ there exists $\varphi_{\epsilon }\in W^{m,2}_0(\Omega_{\epsilon} )$ such that  $\varphi_{\epsilon}\to \varphi $ in $W^{m,2}_0(D )$ as $\epsilon \to 0$.
    \item[{\rm (M2)}] If $v\in W^{m,2}_0(D)$ and there exists a sequence $v_{\epsilon_n}\in W^{m,2}_0(\Omega_{\epsilon_n} )$ such that 
    $v_{\epsilon_n}  \rightharpoonup v  $ in $W^{m,2}(D)$ as $\epsilon_n\to 0 $, then $v\in W^{m,2}_0(\Omega )$.
\end{itemize}
\end{defn}

Note that in the previous definition it is understood that functions $\varphi, \varphi_{\epsilon }$, $v_{\epsilon_n }$  are extended by zero outside their domain of definition. Moreover, the condition  $v  \in W^{m,2}_0(\Omega )$ in (M2) has to be understood in the
sense that the function $v $ can be approximated in $W^{m,2}(D)$ by a sequence of $C^{\infty }$-functions with compact support in $\Omega $. 

We prove the following expected result.

\begin{thm} The Mosco convergence of the spaces $W^{m,2}_0(\Omega_{\epsilon })$ to  $W^{m,2}_0(\Omega )$ as $\eps \to 0$ implies the validity of condition {\bf (C)}, hence the ${\mathcal{E}}$-compact convergence of the operators $H^{-1}_{W^{m,2}_0(\Omega_{\epsilon})}$ to the operator $H^{-1}_{W^{m,2}_0(\Omega )}$.
\end{thm}

{\bf Proof.} Assume that $D$, $\Omega_{\epsilon}$ and $\Omega $ are as in Definition~\ref{Mosco}. We set $K_{\epsilon}=\Omega_{\epsilon }\cap \Omega$.  We divide the proof in several steps. 

{\it Step 1}  We prove that condition
(\ref{limmeas}) hold.  
 Using standard properties of the Lebesgue measure, to prove (\ref{limmeas}) it is enough to show that for any compact set $K\subset \Omega$, we have $|K\setminus \Omega_\eps|\to 0$. But, since $K$ is compact and $K\subset \Omega$, we have the existence of a function $\varphi\in C_0^\infty(\Omega)$ with $\varphi\equiv 1$ in $K$. From (M1) we have a sequence of functions  $\varphi_\eps\in W^{m,2}_0(\Omega_{\eps})$, such that $\varphi_\eps\to \varphi$ in $W^{m,2}(D)$ and in particular in $L^1(D)$.  But this implies that 
$$|\{x\in D:  |\varphi_\eps(x)-\varphi(x)|>1/2\}|\to 0$$
and therefore, since in $K\setminus\Omega_\eps$ we have $\varphi(x)-\varphi_\eps(x)=1$, then $K\setminus\Omega_\eps\subset \{x\in D:  |\varphi_\eps(x)-\varphi(x)|>1/2\}$ and therefore $|K\setminus \Omega_\eps|\to 0$.

{\it Step 2} We prove that condition (C1) is satisfied.  Assume that  $v_{\epsilon }\in W^{m,2}_0(\Omega_{\epsilon })$ are as in (C1) and assume directly that 
$v_{\epsilon }$ are extended by zero outside $\Omega_{\epsilon}$. Since $v_{\epsilon }\in W^{m,2}_0(D )$, by the Sobolev's Embedding Theorem $v_{\epsilon
 }\in L^p(D)$ for some $p>2$. Thus,
$$
\| v_{\epsilon }\|_{L^2(\Omega\setminus K_{\epsilon })}\le |\Omega\setminus K_{\epsilon } |^{\frac{1}{2}-\frac{1}{p}}\| \| v_{\epsilon }\|_{L^p(\Omega\setminus K_{\epsilon })}\le c |\Omega\setminus K_{\epsilon } |^{\frac{1}{2}-\frac{1}{p}}
$$ 
where $c>0$ is a constant independent of $\epsilon$. This, combined with (\ref{limmeas}) implies that $\| v_{\epsilon }\|_{L^2(\Omega \setminus K_{\epsilon })}\to 0$ as $\epsilon \to 0$. Thus, in order to prove that condition (C1) is satisfied, it suffices to prove that   $\| v_{\epsilon }\|_{L^2(\Omega_{\epsilon }\setminus \Omega)}\to 0$ as $\epsilon \to 0$. Assume by contradiction that this is not the case. Then there exists a subsequence 
$ v_{\epsilon_n}$ such that $\| v_{\epsilon_n }\|_{L^2(\Omega_{\epsilon_n }\setminus \Omega)}\to c>0$. Moreover, possibly passing to a subsequence, by (M2)
there exists $v \in W^{m,2}_0(\Omega )$ such that $v_{\epsilon_n}$ converges to $v $ in $L^2(D)$. In particular, it follows that 
$\| v \|_{L^2(D\setminus \Omega )}=c$ which contradicts the fact that $v \in W^{m,2}_0(\Omega)$.
 
{\it Step 3} We prove that the validity of (M1) implies the validity of (C2). For any $\varphi \in W^{m,2}_0(\Omega )$,  we set $T_{\epsilon }\varphi =\varphi_{\epsilon}$ where $\varphi_{\epsilon}$ is as in (M1). Obviously, we have that 
$$
Q_{K_{\epsilon}}(T_{\epsilon }\varphi-\varphi)\le c \| \varphi_{\epsilon}-\varphi  \|_{W^{m,2}_0(D)},
$$
where $c>0$ is independent of $\epsilon$. Since $\varphi_{\epsilon}\to \varphi $ in $W^{m,2}_0(D)$, 
follows that condition (C2) (i) is satisfied. We now prove that condition (C2) (ii) is satisfied. We note that there exists a constant $c>0$ independent of $\epsilon $ such that 
\begin{eqnarray}\lefteqn{
Q_{\Omega_{\epsilon }\setminus K_{\epsilon }}(T_{\epsilon })\le c \| \varphi_{\epsilon } \|_{W^{m,2}(\Omega_{\epsilon }\setminus K_{\epsilon } )} \le c \| \varphi_{\epsilon } -\varphi  \|_{W^{m,2}(\Omega_{\epsilon }\setminus K_{\epsilon } )}  }\nonumber \\
& & \qquad\quad +c \| \varphi \|_{W^{m,2}(\Omega_{\epsilon }\setminus K_{\epsilon } )}  
 = \| \varphi_{\epsilon } -\varphi  \|_{W^{m,2}(\Omega_{\epsilon }\setminus K_{\epsilon } )}+c \| \varphi \|_{W^{m,2}(\Omega \setminus K_{\epsilon } )}\label{mostep3}. 
\end{eqnarray}
It is now clear that by (M1) and the absolute continuity of Lebesgue integrals, the right-hand side of (\ref{mostep3})  goes to zero 
as $\epsilon\to 0$, hence condition (C2) (ii) is satisfied. Condition (C2) (iii) is trivial. 

{\it Step 4} We prove that the validity of (M2) implies the validity of (C3). We set $E_{\epsilon }(v)={\rm Ext}_0v$ where ${\rm Ext}_0$ is the extension-by-zero operator. Clearly, conditions (C3) (i), (ii)  are trivially satisfied. We now consider condition (C3) (iii). Let $v_{\epsilon}$ be as in (C2) (iii). Then 
${\rm Ext}_0v_{\epsilon }$ is a bounded sequence in $W^{m,2}_0(D)$. Thus, there exists $v\in W^{m,2}_0(D)$ such that,  possibly passing to a subsequence, ${\rm Ext}_0 v_{\epsilon }$ is convergent to $v$ strongly in $L^2(D )$ and weakly in $W^{m,2}_0(D )$. By condition (M2), we immediately have that $v\in W^{m,2}_0(\Omega )$, hence condition (C2) (iii) is also satisfied. 

\begin{rem} Sufficient conditions ensuring the Mosco convergence of spaces $W^{m,2}_0(\Omega_{\eps})$ to $W^{m,2}_0(\Omega)$ are well-known. We refer to Bucur and Buttazzo~\cite{bucbut}
and Henrot~\cite{henrot} for a detailed discussion in the case $m=1$. We note that such conditions typically involve geometric notions describing the vicinity of sets (for example, the Hausdorff distance) as well as uniform regularity or topological assumptions on the domains.  Some of the conditions known in the case $m=1$ easily extends to the case $m>1$, as in the case of the compact convergence of sets. For example, if $\Omega_{\eps}\subset \Omega $ is a sequence of open sets compact convergent to $\Omega$ as $\eps \to 0$ (i.e., for any compact set $K\subset \Omega $ there exists $\eps_K>0$ such that $K\subset \Omega_{\eps}$ for all $0<\eps <\eps_K$) then one can prove that the spaces $W^{m,2}_0(\Omega_{\eps})$ converge in the sense of Mosco to $W^{m,2}_0(\Omega)$ as $\eps \to 0$. See Babuska and Vyborny~\cite{babvyb} for more information. 
\end{rem}

\section{Neumann boundary conditions}
\label{neumann}

In this section we consider the operator (\ref{classicop}) subject to Neumann boundary conditions on bounded open sets $\Omega $ in ${\mathbb{R}}^N$.  This has to be understood in the general frame discussed in Section \ref{preliminaries} as follows: by  Neumann boundary conditions we mean that the domain $V(\Omega )$ of the corresponding quadratic form $Q_{\Omega }$ is given by
$$
V(\Omega )=W^{m,2}(\Omega ) ,
$$
and this will be understood throughout this section. Here we assume that the coefficients $A_{\alpha  \beta }$ satisfy the uniform ellipticity condition (\ref{unielp}).

It is well known that both the smoothness of the domains and the kind of perturbations that we are allowed when dealing with operators with Neumann boundary conditions is more restrictive than in the Dirichlet case. An appropriate setting for this issue is clarified with the notion of atlas, as for instance in \cite[Definition 2.4]{bula2012}.  For the sake of completeness and clarity, let us include here the definition.  

For any given set $V\in\mathbb{R}^N$ and $\delta>0$ we denote by $V_\delta$ the set $\{x\in \mathbb{R}^N:  \, d(x,\partial\Omega)>\delta\}$. Moreover, by a cuboid we mean any rotation of a rectangular parallelepiped in $\mathbb{R}^N$.  

\begin{defn}\label{atlas}[Definition 2.4, \cite{bula2012}]
Let $\rho>0$, $s,s'\in \mathbb{N}$ with $s'<s$. Let also $\{V_j\}_{j=1}^s$ be a family of bounded open cuboids and $\{r_j\}_{j=1}^s$ be a family of rotations in $\mathbb{R}^N$.  We say that $\mathcal{A}=(\rho, s,s',\{V_j\}_{j=1}^s,\{r_j\}_{j=1}^s)$ is an atlas in $\mathbb{R}^N$ with parameters $\rho, s,s',\{V_j\}_{j=1}^s,\{r_j\}_{j=1}^s$, briefly an atlas in $\mathbb{R}^N$.  Moreover, we consider the family of all open sets $\Omega\subset \mathbb{R}^N$ satisfying the following:
\par i) $\Omega\subset \cup_{j=1}^s (V_j)_\rho$ and $(V_j)_\rho\cap \Omega\ne \emptyset$ 

\par ii) $V_j\cap\partial\Omega\ne \emptyset$ for $j=1,\dots, s'$ and $V_j\cap\partial\Omega=\emptyset$ for $s'<j\leq s$

\par iii) for $j=1,\dots,s$ we have
$$r_j(V_j)=\{ x\in \mathbb{R}^N:   a_{ij}<x_i<b_{ij}, i=1,\dots, N\}, \quad  j=1,\dots, s$$
$$r_j(V_j\cap \Omega)=\{ x\in \mathbb{R}^N:   a_{Nj}<x_N<g_j(\bar x)\}, \quad j=1,\dots, s'$$
where $\bar x=(x_1,\dots, x_{N-1})$, $W_j=\{ x\in \mathbb{R}^{N-1}:   a_{ij}<x_i<b_{ij}, i=1,\dots, N-1\}$
and the functions $g_j\in C^{k,\gamma}(W_j)$ for $j=1,\dots, s'$, with $k\in\mathbb{N}\cup\{0\}$ and $0\leq \gamma \leq 1$.  Moreover, for $j=1,\ldots, s'$ we have $a_{Nj}+\rho\leq g_j(\bar x)\leq b_{Nj}-\rho$, for all $\bar x\in W_j$. 

We say that an open set $\Omega$ is of class $C^{k,\gamma}_M(\mathcal{A})$ if all the functions $g_j$, $j=1,\dots, s'$ defined above are of class $C^{k,\gamma}(W_j)$ and  $\|g_j\|_{C^{k,\gamma}(W_j)}\leq M$.  We say that an open set $\Omega$ is of class $C^{k,\gamma}(\mathcal{A})$ if it is of class $C^{k,\gamma}_M(\mathcal{A})$ for some $M>0$. Also, we say that an open set $\Omega$ is of class $C^{k,\gamma}$ if it is of class $C^{k,\gamma}_M(\mathcal{A})$ for some atlas $\mathcal{A}$ and some $M>0$. Finally, we  denote by $C^{k}$ the class $C^{k,0}$ for $k\in \mathbb{N}\cup\{0\}$. 

\end{defn}

We recall that if $\Omega $ is a $C^{0}$ bounded open set  then the Sobolev space $W^{m,2}(\Omega )$ is compactly embedded in $L^2(\Omega )$, see e.g., Burenkov~\cite{bu}. Thus, as it is explained in Section 2, the operator $H_{W^{m,2}(\Omega )}$ is well-defined and has compact resolvent.

In this section we discuss the ${\mathcal{E}}$-compact convergence of the operators $H_{W^{m,2}(\Omega_{\epsilon } )}^{-1}$ on families of open sets $\Omega_{\epsilon}$, $\epsilon >0$.
The following theorem is in fact a generalization to higher order operators of \cite{arrJDE95}, \cite[Prop.~2.3]{arca}.

\begin{thm}\label{dumbthm1}
Let $\Omega$ be a bounded  open set in ${\mathbb{R}}^N$ of class $C^{0,1}$ and $\Omega_{\epsilon}$, with $\epsilon >0$,  be bounded open sets in ${\mathbb{R}}^N$ of class $C^{0}$. Let ${\mathcal{A}}$ be an atlas in ${\mathbb{R}}^N$, as in Definition \ref{atlas}, and $M>0$.
Assume that for each $\epsilon >0$ there exists an open set $K_{\epsilon }\subset \Omega \cap \Omega_{\epsilon }$ of class $C^{0,1}_M({\mathcal{A}})$ satisfying (\ref{limmeas}) and one of the following two equivalent conditions:
\begin{itemize}
  \item[i)]  If $v_{\epsilon }\in W^{m,2}(\Omega_{\epsilon })$ and  $\displaystyle\sup_{\epsilon >0}\| v_{\epsilon } \|_{W^{m,2}(\Omega _{\epsilon} )}<\infty$  then  $\displaystyle\lim_{\epsilon \to 0} \| v_{\epsilon }\|_{L^2(\Omega_{\epsilon }\setminus K_{\epsilon })} =0\, ;$
 
  \item[ii)]      $\displaystyle \lim_{\epsilon \to 0}\tau_{\epsilon }=\infty $,  where
$$\displaystyle \tau_{\epsilon }=\inf_{\substack{v\in W^{m,2}(\Omega_{\epsilon })\setminus \{0\} \\ v=0 \ {\rm on}\ K_{\epsilon }}} \frac{Q_{\Omega_{\epsilon }}(v)}{\| v\|^2_{L^2(\Omega_{\epsilon })}} . $$
\end{itemize}
Then we have  $\displaystyle \lim_{\epsilon \to 0}|\Omega_{\epsilon }\setminus K_\eps|=0$  and condition (C) is satisfied. Hence,   $H^{-1}_{W^{m,2}(\Omega_{\epsilon})}\Ccon H^{-1}_{W^{m,2}(\Omega )}$.
\end{thm}

{\bf Proof.}  First we note that conditions i) and ii) both imply that 
$\displaystyle \lim_{\epsilon \to 0}|\Omega_{\epsilon }\setminus K_\eps|=0$.  
The proof of this is very similar as the one from \cite[Prop. 2.3]{arca} and we skip it.

That i) implies ii) is very simple. Indeed, if  $\tau _{\epsilon_k}$ is bounded for some sequence $\epsilon_k\to 0$, then one can find functions $v_{\epsilon_k}\in W^{m,2}(\Omega_{\epsilon_k })$ with $\| v_{\epsilon_k}\|_{L^2(\Omega_{\epsilon_k }\setminus K_{\epsilon_k })}=1$ and $Q_{\Omega _{\epsilon_k}}(v_{\epsilon_k})  $ bounded. Hence, $\|v_{\epsilon_k} \|_{W^{m,2}(\Omega_{\epsilon_k })}$ is a bounded sequence, which contradicts $i)$. 


We now prove that ii) implies i). If i) does not hold, one can find a sequence $\epsilon_k\to 0$ and functions $v_{\epsilon _k}\in  W^{m,2}(\Omega_{\epsilon_k })$ with $\sup_{k\in {\mathbb{N}}}\| v_{\epsilon_k}\|_{W^{m,2}(\Omega_{\epsilon_k })} <\infty $ and $\| v_{\epsilon_k}\|_{L^2(\Omega_{\epsilon_k}\setminus K_{\epsilon_k })}=1$.
Since $K_{\epsilon }$ is of class  $C^{0,1}_M({\mathcal{A}})$ for all $\epsilon >0$, there exists a bounded linear extension operator ${\rm Ext}_{K_{\epsilon } }:W^{m,2}(K_{\epsilon } )\to W^{m,2}({\mathbb{R}}^N )$ with a uniformly bounded norm, {\it i.e.,}
\begin{equation}
\label{uninorm}
\sup_{\epsilon > 0}\bigl\| {\rm Ext}_{K_{\epsilon } }  \bigr\|_{ W^{m,2}(K_{\epsilon } )\to W^{m,2}({\mathbb{R}}^N ) }<\infty, 
\end{equation}
see Burenkov~\cite[Thm.~3,~Chp. 6]{bu} or the classical Stein's book \cite{Stein}.

We set $V_{\epsilon_k}={\rm Ext}_{K_{\epsilon }} ({v_{\epsilon_k}}_{|K_{\epsilon_k}}  ) $. By using the Sobolev's Embedding Theorem and (\ref{uninorm}) one can prove as in \cite{arca} that  $\| V_{\epsilon_k}\|_{L^2(\Omega_{\epsilon_k }\setminus K_{\epsilon_k })}\to 0$.
Consider now the function $w_{\epsilon_k}=V_{\epsilon_k}-v_{\epsilon_k}$. It is clear that $w_{\epsilon_k}=0$ on $K_{\epsilon_k}$,   $\sup_{k\in {\mathbb{N}}}Q_{\Omega_{\epsilon_k}}(w_{\epsilon_k})<\infty$ and $\liminf_{k\to \infty }\|w_{\epsilon_k} \|_{L^2(\Omega_{\epsilon_k} )}>0$. It follows that $\sup_{k\in {\mathbb{N}}}\tau_{\epsilon_k}$ $<\infty $ which contradicts ii).

It remains to prove that condition i) implies that condition (C) is satisfied. By (\ref{unielp}) it follows that the norm $Q_{\Omega }^{1/2}$ is equivalent to the Sobolev norm (\ref{norm}). Thus condition (i) implies the validity of condition {(C1)}.

Since $\Omega$ is of class $C^{0,1}$ there exists a bounded linear extension operator ${\rm Ext}_{\Omega }:W^{m,2}(\Omega )\to W^{m,2}({\mathbb{R}}^N )$. For every $\epsilon >0$, let $T_{\epsilon }$ be the operator of $W^{m,2}(\Omega ) $ to $W^{m,2}(\Omega_{\epsilon } )$ defined by $T_{\epsilon }\varphi = ({\rm Ext}_{\Omega }\varphi )_{|\Omega_{\epsilon }} $, for all $\varphi \in W^{m,2}(\Omega )$. Since $T_{\epsilon }\varphi =\varphi $ on $K_{\epsilon }$
and $T_{\epsilon }\varphi = {\rm Ext}_{\Omega }\varphi $ on $\Omega_{\epsilon }$ for all $\varphi \in W^{m,2}(\Omega )$, and $\lim_{\epsilon \to 0}|\Omega_{\epsilon }\setminus K_{\epsilon}|=0 $, it follows that $T_{\epsilon }$ satisfies condition {(C2)}.

For every $\epsilon >0$, let $E_{\epsilon }$ be the operator of $W^{m,2}(\Omega_{\epsilon} ) $ to $W^{m,2}(\Omega )$ defined by $E_{\epsilon }u= {\rm Ext}_{K_{\epsilon } }(u_{|K_{\epsilon }}) $, for all $u\in W^{m,2}(\Omega_{\epsilon} )$. It is obvious that condition {(C3)} (i) is satisfied. Moreover, by (\ref{uninorm}) also condition {(C3)} (ii) is satisfied. We now prove that { (C3)} (iii) is satisfies as well. Let $v_{\epsilon }\in W^{m,2}(\Omega_{\epsilon })$  and $v\in L^2(\Omega )$ be as in {(C3)} (iii). Since $\| v_{\epsilon }\| _{W^{m,2}(\Omega _{\epsilon })}$
is uniformly bounded, by (\ref{uninorm}) it follows that $ {\rm Ext}_{K_{\epsilon } }(u_{|K_{\epsilon }}) $ is uniformly bounded in $W^{m,2}({\mathbb{R}}^N)$. By the reflexivity of $W^{m,2}({\mathbb{R}}^N)$ there exists $\tilde v\in W^{m,2}({\mathbb{R}}^N)$ such that $v_{\epsilon }$ converges weekly to $\tilde v$ in $W^{m,2}({\mathbb{R}}^N)$ as $\epsilon \to 0$. Thus, $E_{\epsilon }v_{\epsilon } $ converges weekly to $\tilde v$ in $L^{2}(\Omega )$ as $\epsilon \to 0$ hence $\tilde v=v$ and $v\in W^{m,2}(\Omega )$. Thus,  also condition {(C3)} (iii) is satisfied.
The proof is complete. \hfill $\Box$

\par\medskip 
We have the following result, which considers the particular case where $\Omega\subset \Omega_\eps$.  
\begin{corol}\label{auxi}
Let $\Omega$ be a bounded  open set in ${\mathbb{R}}^N$ of class $C^{0,1}$ and $\Omega_{\epsilon}$, with $\epsilon >0$,  be bounded open sets in ${\mathbb{R}}^N$ of class $C^{0}$ with $\Omega\subset \Omega_\eps$.
Assume that one of the following two equivalent conditions is satisfied
\begin{itemize}
\item[i)] If $v_{\epsilon }\in W^{m,2}(\Omega_{\epsilon })$ and  $\sup_{\epsilon >0}\| v_{\epsilon } \|_{W^{m,2}(\Omega _{\epsilon} )}<\infty$  then  $\lim_{\epsilon \to 0} \| v_{\epsilon }\|_{L^2(\Omega_{\epsilon }\setminus \Omega )} =0.$
\item[ii)]  $\displaystyle \lim_{\epsilon \to 0}\ \inf_{\substack{v\in W^{m,2}(\Omega_{\epsilon })\setminus \{0\} \\ u=0 \ {\rm on}\ \Omega }} \frac{Q_{\Omega_{\epsilon }}(v)}{\| v\|^2_{L^2(\Omega_{\epsilon })}}=\infty .$
\end{itemize}

Then we have $\lim_{\epsilon \to 0} |\Omega_{\epsilon }\setminus \Omega |=0
$ and condition (C) is satisfied. Hence,   $H^{-1}_{W^{m,2}(\Omega_{\epsilon})}$ is ${\mathcal{E}}$-compact convergent  to $ H^{-1}_{W^{m,2}(\Omega )}$.
\end{corol}
{\bf Proof.} The proof immediately follows by setting $K_{\epsilon }=\Omega $ for all $\epsilon >0$ and applying Theorem~\ref{dumbthm1}.\hfill $\Box$\\

 We also can prove a simple criterion ensuring spectral stability for Neumann boundary conditions, which generalizes the condition given in  Arrieta and Carvalho~\cite[\S~5.1]{arca}. 
 This criterion can be easily formulated in terms of
the notion of the atlas distance $d_{\mathcal{A}}^{(m)}$ which is introduced in \cite{bulahigh}.   

\begin{defn} {\bf (Atlas distance)}
\label{dev}
Let ${\mathcal{A}} =(\rho , s, s', \{V_j\}_{j=1}^s , \{r_j\}_{j=1}^s )$ be an atlas in ${\mathbb{R}}^N$.
For all $\Omega_1 , \Omega_2\in C^m({\mathcal{A}})$ we set 
\begin{equation}\label{dev1}
d^{(m)}_{{\mathcal{A}}}(\Omega_1,\Omega_2)=\max_{j=1, \dots ,s}\sup_{0\le |\alpha | \le m}\sup_{(\bar x , x_N)\in r_j(V_j)}\left| D^{\alpha }g_{1j}(\bar x) -D^{\alpha } g_{2j}(\bar x)  \right|,
\end{equation}
where $g_{1j}$, $g_{2j}$ respectively, are the functions describing the boundaries of $\Omega_1, \Omega_2$ respectively, as in Definition \ref{atlas}.

Moreover,  we set $d_{{\mathcal{A}}}=d^{(0)}_{{\mathcal{A}}}$ and we  called $d_{{\mathcal{A}}}$  `atlas distance'.  
\end{defn}

The atlas distance clearly depends on the atlas but has the advantage of being easily computable. In the case of open sets of class $C^{0,1}_M$, $d_{\mathcal{A}}$ is equivalent to the usual Hausdorff distance.  Amongst the basic properties of $d_{\mathcal{A}}$,  it is also worth mentioning that $(C(\mathcal{A}), d_{\mathcal{A}})$ is a complete metric space. See \cite{bulahigh}
for more information. 

Then we can prove the following 

\begin{thm}\label{corpoin}
Let ${\mathcal{A}}$ be an atlas in ${\mathbb{R}}^N$  and $\Omega$ be a bounded open set of class $C^{0,1}({\mathcal{A}})$. Let
$\Omega_{\eps}$, $\eps >0,$ be  bounded open sets of class $C({\mathcal{A}})$ such that 
$$ \lim_{\eps \to 0}d_{\mathcal{A}}(\Omega_{\eps}, \Omega)=0.$$
Then condition (C) is satisfied, hence   $H^{-1}_{W^{m,2}(\Omega_{\epsilon})}\Ccon  H^{-1}_{W^{m,2}(\Omega )}$. 
\end{thm}

{\bf Proof.} Let $K_{\eps}$ be bounded open sets of class $C^{0,1}_M({\mathcal{A}})$ with $M>0$ independent of $\eps$, such that $K_{\eps}\subset \Omega\cap \Omega_{\eps}$ and such that $d_{\mathcal{A}}(K_{\eps}, \Omega_{\eps})\le 2 d_{\mathcal{A}}(\Omega _{\eps}, \Omega)$. We denote by $g_{\Omega_{\eps} , j }$ and
$g_{K_{\eps} , j }$ the functions describing the boundaries of $\Omega_{\eps}$  and $K_{\eps}$ respectively, as in Definition~\ref{atlas}. 

 Recall now that if a function $f$ belongs to a Sobolev space the type $W^{m,2}(a,b)$ where $(a,b)$ is a bounded real interval and $f^{(i)}(a)=0$ for any $i=0, \dots , m-1$ then the following Poincar\'{e} inequality holds
\begin{equation}\label{highpoin}
\| f\|_{L^2(a,b)}\le C(b-a)^m\| f^{(m)}\|_{L^2(a,b)},
\end{equation}
where $C>0$ depends only on $m$. 

 Let $v\in W^{m,2}(\Omega_{\eps})$ be such that ${v}_{|K_{\eps}}=0$. Then by applying Fubini-Tonelli's Theorem, using (\ref{highpoin}) and the notation from Definition~\ref{atlas},  we have
\begin{eqnarray}\| v\|^2_{L^2(\Omega_{\eps}\setminus K_{\eps} )}& \le & \sum_{j=1}^{s'}\int_{W_{j}}\int_{g_{K_{\eps} , j }(\bar x) }^{g_{\Omega_{\eps} , j }(\bar x) }|v(\bar x, x_N)|^2dx_Nd\bar x\nonumber \\
&\le  & c  \sum_{j=1}^{s'}\int_{W_{j}}  | g_{\Omega_{\eps} , j }(\bar x)-g_{K_{\eps} , j }(\bar x)  |^{2m}  \int_{g_{K_{\eps} , j }(\bar x) }^{g_{\Omega_{\eps} , j } (\bar x)} \left| \frac{\partial^m v}{\partial x_N^m}    \right|^2dx_nd\bar x  \nonumber \\
& \le &  c d_{\mathcal{A}}^{2m}(K_{\eps}, \Omega_{\eps}) \| v\|^2_{W^{m,2}(\Omega_{\eps})},
\end{eqnarray}
by which we immediately deduce that $\lim_{\eps\to 0}\tau_{\eps}=\infty $, where $\tau_{\eps} $ is defined in Theorem~\ref{dumbthm1}. Thus, Theorem~\ref{dumbthm1} allows to conclude the proof. \hfill $\Box$

\begin{rem}
We note that in Theorem~\ref{corpoin} the assumptions on the open sets $\Omega_{\eps}$ are quite weak. Indeed, it is not required that the sets $\Omega_{\epsilon}$ belong to 
a uniform Lipschitz class and it is only required that they are of class $C({\mathcal{A}})$. In particular, the modulus of continuity of the functions describing their boundaries may 
blow up as $\eps\to 0$. 
\end{rem}

\section{Intermediate boundary conditions }
\label{intermediate}


In this section we consider the operator (\ref{classicop}) subject to intermediate boundary conditions on bounded open sets $\Omega $ in ${\mathbb{R}}^N$ with smooth boundaries. By intermediate boundary conditions we mean that the domain $V(\Omega )$ of the corresponding quadratic form $Q_{\Omega }$ is given by
\begin{equation}\label{quaddom}
V(\Omega )=W^{m,2}(\Omega )\cap W^{k,2}_0(\Omega ),
\end{equation}
with $m\geq 2$, where $k\in {\mathbb{N}}$, $1\le k<m$, is fixed. This will be understood throughout this section. By well-known estimates for intermediate derivatives
(see Burenkov~\cite[p.~160]{bu}) it follows that $V(\Omega )$ is a closed subspace of $W^{m,2}(\Omega )$.

Here we assume that the coefficients $A_{\alpha  \beta }$ are fixed and satisfy the uniform ellipticity condition
(\ref{unielp}).  Thus the operator $H_{V(\Omega )}$ is well-defined and has compact resolvent since   $V(\Omega )$ is compactly embedded in $L^2(\Omega )$.

In this section we discuss the ${\mathcal{E}}$-convergence of the operator $H_{V(\Omega_{\epsilon }) }^{-1}$ on suitable families of smooth open sets $\Omega_{\epsilon }$, $\epsilon >0$. Our analysis includes the case of open sets $\Omega_{\epsilon}$ with oscillating boundaries.

We consider, as in Definition \ref{atlas}, a fixed atlas $\mathcal{A}$ and constant $M>0$ and assume $\Omega\in C^{m}_M(\mathcal{A})$.  We will also consider that $\Omega_\eps\in C^{m}_{M_\eps}(\mathcal{A})$ for some constants $M_\eps$ not necessarily uniformly bounded in $\epsilon$. To simplify the proofs of the results, we will consider that the perturbation of the boundary is localized in just one of the cuboids $V_i$ for some $i=1,\ldots,s'$, that is, one of the cuboids which touch the boundary.  We refer to Corollary \ref{weaklemint2}  for a general statement when the perturbation acts not in just a single cuboid.

Therefore, let us denote the cuboid by $V$ (we will drop the subindex)  and hence we will assume that $\Omega\setminus V_\rho=\Omega_\eps\setminus V_\rho$, that is, the perturbation is localized in the interior of $V$.  Without loss of generality, we may assume that 
$V=W\times (a,b)$, where $W=\{x\in \mathbb{R}^{N-1}:  a_j<x_j<b_j, j=1,\ldots,N-1\}$. Moreover,
$V\cap \Omega= \{ (\bar{x},x_N)\in W\times ]a,b[ \; : \;    a <x_N<g(\bar{x}) \}$ and
$V\cap \Omega_{\epsilon } = \{ (\bar{x},x_N)\in  W\times ]a,b[ \; : \;   a <x_N<g_{\epsilon }(\bar{x}) \}$. 
The functions $g$ and $g_\eps$ define the boundary of $\Omega$ and $\Omega_\eps$ in $V$ and  as in Definition~\ref{atlas} we assume that $a+\rho <g, g_{\eps} <b-\rho $,     $g,g_\eps \in C^m(W)$ with $\|g\|_{C^m(W)}\leq M$,  $\|g_\eps\|_{C^m(W)}\leq M_\eps$.

The following lemma provides a sufficient condition for the ${\mathcal{E}}$-compact convergence of  $H^{-1}_{V(\Omega_{\epsilon})}$ to $H^{-1}_{V(\Omega)}$.  As usual, $\| \cdot \|_{\infty}$ denotes the $L^{\infty }$-norm.

\begin{lem}\label{weaklemint} With the notation above and assuming that for every $\epsilon >0$ there exists $\kappa_{\epsilon }>0$ such that 
\begin{itemize}
\item[(i)] $\kappa_{\epsilon }> \| g_{\epsilon }-g\|_{\infty},\ \ \forall\ \epsilon >0$ 
\item[(ii)] $\lim_{\epsilon\to 0}\kappa_{\epsilon }=0$
\item[(iii)]  $\lim_{\epsilon \to 0}\frac{\| D^{\beta }(g_{\epsilon }-g)\|_{\infty }}{\kappa_{\epsilon }^{m-|\beta |-\frac{1}{2}}}=0,$
 $\forall \ \beta \in {\mathbb {N}}^N$ with $|\beta|\le m$.
\end{itemize}

Then condition {\bf (C)} is satisfied, hence  $H^{-1}_{V(\Omega_{\epsilon})}\Ccon H^{-1}_{V(\Omega )}$.
\end{lem}

{\bf Proof.}   
 Our argument is based on the construction of a suitable diffeomorphism from $\bar\Omega_{\epsilon }$ onto $\bar\Omega $ which coincides with the identity outside $V_{\rho}$. Thus, since $\Omega_{\epsilon }\setminus V_{\rho }=\Omega \setminus V_{\rho }$ for all $\epsilon >0$,  we can carry out our construction in $V$  and assume directly that $\Omega=\Omega\cap V$ and $\Omega_\eps=\Omega_\eps \cap V$. Hence, $\Omega = \{ (\bar{x},x_N)\in W\times ]a,b[ \; : \;    a <x_N<g(\bar{x}) \}$ and
$\Omega_{\epsilon } = \{ (\bar{x},x_N)\in  W\times ]a,b[ \; : \;   a <x_N<g_{\epsilon }(\bar{x}) \}$ where $W$ is defined above.  Let\footnote{the value of $ \hat k$ does not play any significant role and is used only to prove (\ref{leint0})} $   \hat k >2(m+1)$. We set $k_{\epsilon }=\hat k\kappa_{\epsilon}$ and  $\tilde g_{\epsilon }= g_{\epsilon}-k_{\epsilon }$, $K_{\epsilon }=  \{ (\bar{x},x_N)\in  W\times ]a,b[ \; : \;   a <x_N<\tilde g_{\epsilon }(\bar{x}) \}$. 
Note that by (i), (ii), $\emptyset \ne K_{\epsilon } \subset \Omega\cap \Omega _{\epsilon }$ for all $\epsilon >0 $ sufficiently small. Moreover, by (ii) condition (\ref{limmeas}) is satisfied.

 Throughout the proof, we will denote by $C$ a generic constant which will be independent of $\eps$ and all the functions involved. This constant may change from line to line. If at some point we want to distinguish some constant we will use another notation and make this clear.  

We now prove that condition (C1) is satisfied. Let $v_{\epsilon }$, $\epsilon >0$, be as in condition (C1). We will need the one dimensional embedding estimates   $\|f\|_{L^\infty(a,b)}\leq K\|f\|_{W^{1,2}(a,b)}$ where the constant $K=K(d)$ is uniformly bounded for $|b-a|\geq d$ (see, e.g., Burenkov~\cite{bu}). Hence, by Tonelli's Theorem, and applying this last estimate to the function $v_\eps(\bar x, \cdot)$ in the interval $(a, g_\eps(\bar x))$,  we get
\begin{equation}\label{1dim}
\begin{array}{l}
\displaystyle \int_{\Omega_{\epsilon }\setminus K_{\epsilon }}|v_{\epsilon }|^2dx=\int_{W}\int_{\tilde g_\eps(\bar x)}^{g_\eps(\bar x)}|v_\eps(\bar x, x_N)|^2dx_n d\bar x \\  \\ 
\displaystyle\qquad \le C\int_W|\tilde g_\eps(\bar x)-g_\eps(\bar x)|  \|v_\eps(\bar x, \cdot) \|_{W^{1,2}(a, g_\eps(\bar x))}^2 d\bar x
\\ \\
\displaystyle\qquad\le C\|\tilde g_\eps-g_\eps\|_{L^\infty(W)}\| v_\eps\|^2_{W^{1,2}(\Omega_\eps)}\leq C \kappa_\eps^2 \ Q_{\Omega_{\epsilon }}(v_{\epsilon })\eto 0.
\end{array} 
\end{equation}

By (\ref{1dim}) the validity of condition (C1) follows.

We now prove that condition (C2) is satisfied.   Let $\Phi_\eps: \bar \Omega_{\epsilon }\to \bar \Omega$  be the map defined by
$\Phi_{\epsilon }(\bar x, x_N)=(\bar x, x_N-h_{\epsilon }(\bar x , x_N))$ for all $(\bar x, x_N)\in \bar\Omega_{\epsilon }$ where
\begin{equation}
\label{acca}
h_{\epsilon }(\bar x, x_N)=\left\{\begin{array}{ll}0,\ \ &{\rm if }\ a\le x_N\le \tilde g_{\epsilon }(\bar x), \vspace{2mm}  \\
(g_{\epsilon }(\bar x)-g(\bar x))\left( \frac{x_N-\tilde g_{\epsilon }(\bar x)}{g_{\epsilon }(\bar x)-\tilde g_{\epsilon }(\bar x)}   \right)^{m+1},\ \ &{\rm if }\  \tilde g_{\epsilon }(\bar x)<x_N\le  g_{\epsilon }(\bar x)\, .
\end{array} \right.
\end{equation}
Geometrically speaking, the map $\Phi_\eps$ transforms the ``vertical segment'' $l_{\bar x}^\eps=\{(\bar x, x_N): a <x_N<g_\eps(\bar x)\}$ to the segment $l_{\bar x}=\{(\bar x, x_N): a<x_N<g(\bar x)\}$, leaving invariant the part of the segment with $a<x_N<\tilde g_\eps(\bar x)$.

Note that   for $\eps$ fixed the transformation $\Phi_{\epsilon }$ is a diffeomorphism of class $C^m$ from $\bar\Omega_{\epsilon }$ onto $\bar\Omega $.  The $C^m$ norm of $\Phi_\eps$ will not be bounded in general for $m\geq 2$ as $\eps\to 0$, but  by the choice of $\hat k$  and having in mind (i) (ii) and (iii), it follows that 

\begin{equation}\label{leint0}
C^{-1}\le |{\rm det} D\Phi_{\epsilon }|\le C.
\end{equation}
 To see this, just note that $|{\rm det} D\Phi_{\epsilon }|=|1-\partial_{x_N} h_\eps|$.  
 
 In order to estimate the derivatives (up to order $m$) of the transformation $\Phi_\eps$, we need to study the derivatives of the function $h_\eps$. By the Leibniz formula we have
$$
D^{\alpha }h_{\epsilon }(x)=\sum_{0\le \gamma \le \alpha }{\alpha \choose \gamma } D^{\gamma }(g_{\epsilon }(\bar x)-g(\bar x))D^{\alpha -\gamma }  \left(\frac{x_N-\tilde g_{\epsilon }(\bar x)}{g_{\epsilon }(\bar x)-\tilde g_{\epsilon }(\bar x)} \right)^{m+1},
$$
for all $(\bar x ,x_N)\in W\times ]a,b[$ with $\tilde g_{\epsilon }(\bar x)<x_N<  g_{\epsilon }(\bar x)$.
By standard calculus, it is easy to check that  
\begin{eqnarray}\left|
D^{\alpha -\gamma }  \left(\frac{x_N-\tilde g_{\epsilon }(\bar x)}{g_{\epsilon }(\bar x)-\tilde g_{\epsilon }(\bar x)} \right)^{m+1}\right|\le \frac{C}{|g_{\epsilon }(\bar x)-\tilde g_{\epsilon }(\bar x) |^{|\alpha |-|\gamma |}}
\leq \frac{C}{\kappa_{\epsilon }^{|\alpha |-|\gamma |}},
\end{eqnarray}
hence,  
\begin{equation}\label{chain0}
\|   D^{\alpha }h_{\epsilon}\|_{\infty }\le C\sum_{0\le \gamma\le \alpha }\ 
\frac{\| D^{\gamma }(g_{\epsilon }-g)\|_{\infty }}{\kappa_{\epsilon }^{|\alpha |-|\gamma |}}
\end{equation}
for all $\epsilon >0$ sufficiently small. 

 Let $T_{\epsilon }$ be the map from $V(\Omega )$
to $V(\Omega_{\epsilon })$ defined by
$$
T_{\epsilon }\varphi =\varphi \circ \Phi_{\epsilon },
$$
for all $\varphi \in V(\Omega )$. Note that $T_{\epsilon }$ is well-defined since $\Phi_{\epsilon }$ is a diffeomorphism of class $C^m$.

Condition (C2) (i) is immediately satisfied since $T_{\epsilon }\varphi =\varphi $ on $K_{\epsilon }$ for all $\varphi \in V(\Omega )$.
We now prove that condition (C2) (ii) is satisfied. Let $\varphi \in V(\Omega )$.
By the chain rule, for any multindex $\alpha $ with $|\alpha| = m$, we have
\begin{equation}
\label{chain1}
D^{\alpha}(\varphi( \Phi_{\epsilon }(x)))=\sum_{1\le |\beta |\le m}D^{\beta }\varphi ( \Phi_{\epsilon}(x)) p_{m,\beta }^{\alpha }(\Phi_{\epsilon} )
\end{equation}
where $ p_{m,\beta }^{\alpha}  (\Phi_{\epsilon} )$ is a homogeneous polynomial of degree $|\beta |$ in derivatives of $\Phi_{\epsilon }$ of order not exceeding $m-|\beta |+1$, and coefficients depending on $\alpha $ but  not depending on $\epsilon$.  
 Since $\Phi_\eps(x)=x-(0,h_\eps(x))$, then a derivative of $\Phi_\eps$ of order $|\beta|$ is either constantly 1 or 0 or a derivative or $h_\eps$ of order les or equal than $|\beta|$. Then $p_{m,\beta }^{\alpha}(\Phi_{\epsilon} )$ is a polynomial of degree less or equal than $|\beta|$ in the derivatives of $h_\eps$ of order not exceeding $m-|\beta|+1$.

In particular, using (\ref{chain0}) this implies that 
$$
\| p_{m,\beta}^{\alpha}(\Phi_\eps)\|_{\infty}\leq C\left( \sum_{0\le | \gamma | \le m-|\beta|+1 }
\frac{\| D^{\gamma }(g_{\epsilon }-g)\|_{\infty }}{\kappa_{\epsilon }^{m-|\beta|+1-|\gamma |}}+1\right)^{|\beta|}$$

Note that,  in particular if $|\beta |=1$ then  
\begin{equation}\label{beta=1}
\| p_{m,\beta}^{\alpha} (\Phi_\eps)\|_{\infty}\leq C\left( \sum_{0\le | \gamma | \le m }
\frac{\| D^{\gamma }(g_{\epsilon }-g)\|_{\infty }}{\kappa_{\epsilon }^{m-|\gamma |}}+1\right)=o(1)\kappa_\eps^{-1/2},
\end{equation}
where, as it is customary, we denote by $o(1)$ a function which goes to 0 as $\eps\to 0$.   

If $2\leq |\beta|\leq m$, we have 
\begin{equation}\label{beta>1}
\| p_{m,\beta}^{\alpha}(\Phi_\eps)\|_{\infty}\leq C\left( \sum_{0\le | \gamma | \le m-1}
\frac{\| D^{\gamma }(g_{\epsilon }-g)\|_{\infty }}{\kappa_{\epsilon }^{m-1-|\gamma |}}+1\right)^{m}\leq C
\end{equation}
where we have used hypothesis (iii).

%

Hence, we have 
\begin{eqnarray}\lefteqn{
Q_{\Omega_{\epsilon }\setminus K_{\epsilon }}(T_{\epsilon}(\varphi ) )
\le   \int_{\Omega_{\epsilon}\setminus K_{\epsilon}} |\varphi (\Phi_{\epsilon})|^2dx +
C
\sum_{|\alpha |=m}\int_{\Omega_{\epsilon}\setminus K_{\epsilon}}
|D^{\alpha }(\varphi (\Phi_{\epsilon}))|^2dx} \nonumber   
\\ & & \le   C\int_{\Omega \setminus K_{\epsilon}} |\varphi |^2dx+
C\sum_{   \substack{|\alpha |=m\\ 1\le |\beta| \le m}      }\int_{\Omega_{\epsilon}\setminus K_{\epsilon}}
|  D^{\beta } \varphi ( \Phi_{\epsilon}(x))   p_{m,\beta }^{\alpha}  (\Phi_{\epsilon} )    |^2dx \nonumber
\\  & & \le   C\int_{\Omega \setminus K_{\epsilon}} |\varphi |^2dx+
C\sum_{   \substack{|\alpha |=m\\ 1\le |\beta| \le m}    } \|    p_{m,\beta }^{\alpha}      (\Phi_{\epsilon} )\|_{\infty}^2\int_{\Omega_{\epsilon}\setminus K_{\epsilon}}
|  D^{\beta } \varphi ( \Phi_{\epsilon}(x))  |^2dx 
\nonumber
\\  & & \label{last-equation}
\le   C\int_{\Omega \setminus K_{\epsilon}} |\varphi |^2dx+
C\sum_{ \substack{|\alpha |=m\\ 1\le |\beta| \le m}    } \|   p_{m,\beta }^{\alpha}    (\Phi_{\epsilon} )\|_{\infty}^2
\ \|\varphi\|_{W^{|\beta|,2}(\Omega\setminus K_\eps)}^2.
\end{eqnarray}

Notice that since $\varphi\in W^{m,2}(\Omega)$ is a fixed function and $|\Omega\setminus K_\eps|\to 0$, then  
\begin{equation}\label{l2to0}
\displaystyle \|\varphi\|_{W^{m,2}(\Omega\setminus K_\eps)}\eto 0 \end{equation}
Also, notice that since $\varphi\in W^{m,2}(\Omega)\hookrightarrow W^{2,2}(\Omega)$, we have that for all $i=1,\ldots,N$, $\nabla \varphi\in W^{1,2}(\Omega)$. Therefore, $\nabla\varphi(\bar x, \cdot)\in W^{1,2}(a, g(\bar x))\hookrightarrow L^\infty(a, g(\bar x))$ a.e. $\bar x\in W$.  With a similar argument as the one we use to show (\ref{1dim}), we have 
\begin{equation}\label{H1-to 0}
\begin{array}{l}
\displaystyle\|\varphi\|^2_{W^{1,2}(\Omega\setminus K_\eps)}\leq C\|\tilde g_\eps-g\|_\infty \|\varphi\|_{W^{2,2}(\Omega)}^2\leq C\kappa_\eps .
%
%
\end{array}
\end{equation}

Hence, from (\ref{l2to0}) we have that the first term in (\ref{last-equation}) goes to 0.  Moreover, for the second term in (\ref{last-equation}) we consider the sum for $|\beta|=1$ and $2<|\beta|\leq m$ separated, apply (\ref{beta=1}) and (\ref{beta>1}) to obtain
$$
\begin{array}{l}
\displaystyle \sum_{  \substack{|\alpha |=m\\ 1\le |\beta| \le m}    } \|  p_{m,\beta }^{\alpha}  (\Phi_{\epsilon} )\|_{\infty}^2
\ \|\varphi\|_{W^{|\beta|,2}(\Omega\setminus K_\eps)}^2 \leq o(1)\kappa_\eps^{-1}\|\varphi\|_{W^{1,2}(\Omega\setminus K_\eps)}^2+ C\|\varphi\|_{W^{m,2}(\Omega\setminus K_\eps)}^2 \eto 0
\end{array}
$$
where we use (\ref{l2to0}) and (\ref{H1-to 0}). 
This shows that (C2) (ii) is satisfied.

Condition (C2) (iii) is trivial. 

We now prove that condition (C3) is satisfied.  By condition (iii) it follows that $\|\nabla g_{\eps}  \|_{\infty }$ is uniformly bounded for $\eps $ sufficiently small (recall that $m\geq 2$).  Thus the open sets $\Omega_{\epsilon }$ belong to the same class  $C^{1}_M(\mathcal{A})$ for a suitable fixed $M>0$.  Hence, there exists a bounded linear  extension operator ${\rm Ext }_{\Omega_{\epsilon} }$ from $W^{m,2}(\Omega_{\epsilon } )$ to $W^{m,2}(\mathbb{R}^N)$ such that  \begin{equation}
\label{uninormbis}
\sup_{\epsilon > 0}\bigl\| {\rm Ext}_{\Omega_{\epsilon } }  \bigr\|_{ W^{m,2}(\Omega_{\epsilon } )\to W^{m,2}( {\mathbb{R}}^N ) }<\infty .
\end{equation}
We set $E_{\epsilon }u=({\rm Ext}_{\Omega_{\epsilon } }u)_{|\Omega }$
for all $u\in V(\Omega_{\epsilon })$. It is straightforward that $E_{\epsilon }$ satisfies conditions (C3) (i), (ii). We now prove that condition (C3) (iii) is satisfied. Let $v_{\epsilon }$, $\epsilon >0$, and $v$ be as in condition (C3) (iii).
Since $E_{\epsilon }v_{\epsilon }$ is bounded
in $W^{m,2}(\Omega ) $ and $W^{m,2}(\Omega )$ is compactly embedded in $L^2(\Omega )$, it follows that there exists $\tilde v\in W^{m,2}(\Omega )$ such that, by possibly considering a subsequence,  $v_{\epsilon } $ converges weakly to $\tilde v $ in $W^{m,2}(\Omega )$ and strongly in $L^2(\Omega )$ as $\epsilon \to 0$. It follows that $v=\tilde v$, hence $ v\in W^{m,2}(\Omega )$. It remains to prove that $\tilde v\in W^{k,2}_0(\Omega )$. To do so it suffices  use the extension-by-zero operator ${\mathcal{E}}$. Indeed, ${\mathcal{E}}(v_{\epsilon })\in W^{k,2}_0(U)$ where $U$ is a bounded open set containing all sets $\Omega_{\epsilon }$ and $\Omega $.
Since ${\mathcal {E}}v_{\epsilon }$ is bounded
in $W^{k,2}_0(U ) $ and $W^{k,2}_0(U )$ is compactly embedded in $L^2(U )$ it follows that there exists $\hat v\in W^{k,2}_0(\Omega )$ such that by possibly considering a subsequence  $v_{\epsilon } $ converges weakly to $\hat v $ in $W^{k,2}_0(U )$ and strongly in $L^2(U)$ as $\epsilon \to 0$. Clearly,   $\hat v$ vanishes outside $\Omega $, hence $\hat v\in W^{k,2}_0(\Omega )$. Obviously, $\hat v =v$ in $\Omega$. The proof is complete. \hfill $\Box$\\

\begin{rem}\label{general-cuboid}
It is not difficult to see that the fact that the set $W$ is an $(N-1)$-dimensional cuboid of the form $W=\{x\in \mathbb{R}^{N-1}:  a_j<x_j<b_j, j=1,\ldots,N-1\}$  is not essential at all in the proof of Lemma \ref{weaklemint}.  As a matter of fact, exactly the same proof works if we consider a general smooth, say piecewise $C^1$,  set $W$ and $V$ the cylinder of base $W$ that is $V=W\times (a,b)$. 
\end{rem}

\begin{rem}\label{gagliardo}
By the classical Gagliardo-Nirenberg  interpolation inequality $\| D^{\beta}f  \|_{\infty}\le C( \sum_{|\alpha |=m}\| D^{\alpha }f\|_{\infty})^{|\beta|/m} \| f \|_{\infty }^{1-|\beta|/m}$  (cf. e.g., \cite[p.~125]{nire}), it turns out that in order to verify condition (iii) in Lemma~\ref{weaklemint}, it suffices to verify it for
$|\beta |=0$ and $|\beta |=m$.
\end{rem}

We can deduce now the following, 
\begin{prop}
With the notations above if $\|g_\eps-g\|_\infty\eto 0$
and if $$\sup_{|\alpha|=m} \{  \|g_\eps-g\|_\infty^{\frac{1}{2m-1}} \,\, 
\|D^\alpha g_\eps- D^\alpha g\|_\infty \}\eto 0,$$
then condition {\bf (C)} is satisfied.  Hence  $H^{-1}_{V(\Omega_{\epsilon})}\Ccon H^{-1}_{V(\Omega )}$.
\end{prop}

{\bf Proof.}  
 Let us denote by $\delta_\eps>0$ a sequence such that 
$\sup_{|\alpha|=m} \{  \|g_\eps-g\|_\infty^{\frac{1}{2m-1}} \,\, 
\|D^\alpha g_\eps- D^\alpha g\|_\infty \}\leq \delta_\eps \eto 0$
and define 
$$\rho_\eps=\max\{ \delta_\eps, \|g_\eps-g\|_{\infty}^{\frac{1}{2m-1}}\}\eto 0$$

Let us choose 
$$\kappa_\eps=\frac{\|g_{\eps}-g\|_{\infty}^{\frac{2}{2m-1}}}{\rho_\eps}$$
and assume directly that $\kappa_\eps\ne 0$. 
Then,  since $2/(2m-1)<1$,  $\|g_\eps-g\|_\infty \to 0$ and $\rho_\eps\eto 0$, for $\eps$ small enough we have
$$\|g_\eps-g\|_\infty<\|g_\eps-g\|_\infty^\frac{2}{2m-1} <\kappa_\eps \leq \|g_\eps-g\|_\infty^\frac{1}{2m-1}\eto 0$$
and therefore hypothesis i) and ii) from  Lemma~\ref{weaklemint} hold.  

Moreover, with the definition of $\kappa_\eps$, $\delta_\eps$ and $\rho_\eps$ and noting that $\rho_\eps\geq \delta_\eps$,  we have 
$$\sup_{|\alpha|=m} \{  \kappa_\eps^{1/2} \,\, 
\|D^\alpha g_\eps- D^\alpha g\|_\infty \}\leq \rho_\eps^{-1/2}\delta_\eps \leq \delta_\eps^{1/2}\eto 0 $$
Using now Remark \ref{gagliardo} we can easily show that iii) from Lemma~\ref{weaklemint} holds and the proposition follows.  \hfill $\Box$\\

%
%
%
%
%

In case the perturbation does not act only over one cuboid or it is different from the set $V$ considered in Remark \ref{general-cuboid} above, we can also prove the following result with the aid of a partition of unity subordinated to the family of cuboids $\{V_i\}_{i=1}^s$.  We will denote the functions that define the boundary of $\Omega$ and $\Omega_\eps$ in $V_i$ for all $i=1,\ldots,s'$ by $g_i $ and $g_{\eps,i}$ respectively.

\begin{corol}\label{weaklemint2} 

With the notations above,  if we assume that for every $\epsilon >0$ there exists $\kappa_{\epsilon }>0$ such that 
\begin{itemize}
\item[(i)] $\kappa_{\epsilon }\geq \| g_{\epsilon, i}-g_i\|_{\infty},\ \ \forall\ \epsilon >0, \quad \forall i=1,\ldots, s'.$ 
\item[(ii)] $\lim_{\epsilon\to 0}\kappa_{\epsilon }=0$
\item[(iii)]  $\lim_{\epsilon \to 0}\frac{\| D^{\beta }(g_{\epsilon,i}-g_i)\|_{\infty }}{\kappa_{\epsilon }^{m-|\beta |-\frac{1}{2}}}=0,$
 $\forall \ \beta \in {\mathbb {N}}^N$ with $|\beta|\le m$ and for all $i=1,\ldots,s'$. 
\end{itemize}

Then condition {\bf (C)} is satisfied, hence  $H^{-1}_{V(\Omega_{\epsilon})}\Ccon H^{-1}_{V(\Omega )}$.
\end{corol}


\begin{rem}\label{gagliardo2}
Observe that a similar observation as in Remark \ref{gagliardo} can be applied in this case. 
\end{rem}

Finally, we can deduce the following
\begin{thm}Let ${\mathcal{A}}$ be an atlas in ${\mathbb{R}}^N$, $M>0$, $m\in {\mathbb{N}}$, $m\geq 2$. Let $\Omega_{\eps}, \Omega\in C^m_M(\mathcal{A})$, $\eps >0$, be such that 
$$
\lim_{\eps\to 0}d_{\mathcal{A}}^{(m-1)}(\Omega_{\eps}, \Omega )=0.
$$
Then condition {\bf (C)} is satisfied, hence  $H^{-1}_{V(\Omega_{\epsilon})}\Ccon H^{-1}_{V(\Omega )}$.
\end{thm}

{\bf Proof.} Apply Corollary~\ref{weaklemint2} with $\kappa_{\eps}=(d_{\mathcal{A}}^{(m-1)}(\Omega_{\eps}, \Omega ))^{\frac{1}{m}}$. \hfill $\Box$

%
%
%
%
%
%
%
%
%
%

\section{The biharmonic operator with intermediate boundary conditions}
\label{biharmonic-intermediate}

In this section, we shall consider the biharmonic operator subject to intemediate boundary conditions in a family of domains with oscillatting boundaries.

 Without loss of generality and to simplify the exposition, let us assume that our domain $\Omega\subset \mathbb{R}^N$ is of the form $\Omega =W\times (-1,0)$ where $W\subset \mathbb{R}^{N-1}$ is either an $(N-1)$ dimensional cuboid as in Lemma \ref{weaklemint} or a smooth domain as in Remark \ref{general-cuboid}.    We also assume  that the perturbed domain $\Omega_{\eps}$ is given by  $\Omega_\eps=\{(\bar x, x_N):   \bar x\in W, \,\, -1<x_N <g_\eps(\bar x)\}$ where  $g_\eps(\bar x)=\eps^\alpha b(\bar x/\eps)$ for all $\bar x \in W$ and $b:\mathbb{R}^{N-1}\to [0, \frac12)$ is a $Y$-periodic fixed function, where  $Y$ is the unit cell  $Y=(-\frac12,\frac12)^{N-1}$ and $\alpha>0$.  Note that for simplicity, we have assumed that $b\geq 0$.   Note that the general case can also be treated in a straightforward way, although to avoid annoying technicalities and for the sake of the exposition we will stick to the simpler case $b\geq 0$. 
 This implies in particular that $\Omega\subset \Omega_\eps$, $\eps>0$. 
We denote by   $\Gamma_{\epsilon }$ the set $\{(\bar x, x_N):\ \bar x \in W, \  x_N=g_{\eps}(\bar x) \}$ which is the part of the boundary of $\Omega_{\epsilon }$ above $W$. 
It is also covenient to set  $\Omega _0=\Omega$. In the sequel we shall also identify $\Gamma _0=W\times \{0\}$ with $W$.

Namely, we shall consider the operators 
$$
H_{\Omega_{\eps} , I}= \Delta^2 u +u
$$
on the  open sets $\Omega_{\eps}$, 
with $u$  subject to the classical boundary conditions
$$
u=0,\ \ {\rm and}\ \  \Delta u -K \frac{\partial u}{\partial \nu }=0,\ \ {\rm on }\ \ \partial \Omega_{\eps} ,
$$
where  $K$ denotes  the mean curvature of $\partial \Omega_{\eps}$, {\rm i.e.} the sum of the principal curvatures.
More precisely, the operators  $H_{\Omega_{\eps} , I} $  are the operators associated with the quadratic form 
\begin{equation}\label{biquadratic}
Q_{\Omega_{\eps} }(u,v)=\int_{\Omega_{\eps} }(D^2u: D^2v  + uv )dx 
\end{equation}
defined for all $u, v\in W^{2,2}(\Omega_{\eps} )\cap W^{1,2}_0(\Omega_{\eps} )$, as discussed in Section 2. Recall that  $D^2u$ denotes the Hessian
matrix of $u$ and $D^2u: D^2v=\sum_{i,j=1}^N\frac{\partial^2u}{\partial x_i\partial x_j}\frac{\partial^2v}{\partial x_i\partial x_j}$. 

It is clear that this is a special case of those discussed in the previous section with $m=2$ and $k=1$ in (\ref{quaddom}). 

It will be convenient to denote by $H_{\Omega , D}$ the operator $\Delta ^2+I$ subject to Dirichlet boundary conditions on $W$, that  is, 
\begin{equation}\label{dirconw}
u=\frac{\partial u}{\partial \nu }=0,\ \ {\rm on}\ \ W
\end{equation}
and intermediate boundary conditions on $\partial \Omega \setminus W$, 
which is the operator canonically associated with the quadratic form (\ref{biquadratic}) defined for all $u, v\in W^{2,2}_{0, W}(\Omega)$ where $W^{2,2}_{0, W}(\Omega)$ is the space  of functions $u$ in $W^{2,2}(\Omega)\cap W^{1,2}_0(\Omega)$ satisfying conditions (\ref{dirconw}).


\par\medskip

We have the following result.

\begin{thm}\label{main}
With the notations above, we have the following trichotomy:
\begin{itemize}
\item[i)] If $\alpha>3/2$, then $H_{\Omega_\eps, I}^{-1}\Ccon H_{\Omega, I}^{-1}$. 
\item[ii)] If $0<\alpha<3/2$ and $b$ is non-constant, then $H_{\Omega_\eps, I}^{-1}\Ccon H_{\Omega, D}^{-1}$. 
\item [iii)] If $\alpha=3/2$, then $H_{\Omega_\eps, I}^{-1}\Ccon \hat H_\Omega ^{-1}$, where $\hat H_\Omega$ is the  operator $\Delta^2+I$ in $\Omega$ with intermediate boundary conditions on $\partial \Omega \setminus W$ and the following boundary conditions in $W$:  $u=0,\quad \Delta u-\gamma \frac{\partial u}{\partial \nu}=0$, 
 where the factor $\gamma$ is given as 
 \begin{equation}\label{maingamma}\gamma=\int_{Y\times (-\infty , 0)  }|D^2 V|^2dy       = -\int_{Y}b\frac{\partial}{\partial y_N}(\Delta_{\bar y} V+\Delta V) d\bar y .
\end{equation} 
Here $\Delta_{\bar y}$ is the Laplace operator in the $\bar y$-variables and the function $V$ is $Y$-periodic in the variables $\bar y$ and  satisfies the following microscopic problem
 \begin{equation}
\left\{
 \begin{array}{ll}
 \Delta^2V=0,\  & \hbox{\rm in } Y\times (-\infty, 0), \\
 V(\bar y, 0)=b(\bar y),\ & \hbox{\rm on } Y, \\
 \frac{\partial^2 V}{\partial y_N^2}(\bar y,0)=0,\  &\hbox{\rm on } Y.
  \end{array}
\right.
 \end{equation}
\end{itemize}
\end{thm}

\begin{rem}
For completeness, we have stated in Theorem \ref{main} the three different cases but up to now we can only show part $i)$ and $ii)$.  We will provide a proof of these two cases and will leave the proof of case $\alpha=3/2$, which is quite involved,  for Section \ref{critical}. 
\end{rem}

{\bf Proof of Theorem \ref{main} $i)$ and $ii)$.} (i)  Let $\tilde\alpha\in ]3/2, \alpha [$.  It is easily verified that conditions (i), (ii), (iii) in Lemma~\ref{weaklemint}  are satisfied with $\kappa_{\epsilon }=\epsilon ^{2\tilde\alpha /3}$ for $\epsilon $ small enough. 

(ii) We prove that condition (C) is satisfied with  $V(\Omega ) =W^{2,2}_{0,W}(\Omega)$ and $V(\Omega_{\epsilon})= W^{2,2}(\Omega_{\epsilon})\cap W^{1,2}_0(\Omega_{\epsilon})$.  Let $K_{\epsilon}=\Omega$. 
 Notice that condition (\ref{limmeas}) is trivially satisfied. Moreover, it is easy to see that (C1) is also satisfied. Since $V(\Omega_\eps)$ is continuously embedded into $W^{1,2}_0(\Omega_\eps)$, we have that if $v_\eps\in V(\Omega_\eps)$ with $Q_{\Omega_\eps}(v_\eps)\leq C$ with $C$ independent of $\epsilon$, then $\|v_\eps\|_{W^{1,2}(\Omega_\eps)}\leq \tilde C$ for some  $\tilde C$ independent of $\epsilon$ (see also Burenkov~\cite[Thm.~6,~ p.~160]{bu}). Using Poincar\'e inequality in the $x_N$ direction in $\Omega_\eps\setminus \Omega$, we easily get $\|v_\eps\|_{L^2(\Omega_\eps\setminus\Omega)}^2\leq \rho(\eps)\|\partial_{x_N}v_\eps\|_{L^2(\Omega_\eps\setminus\Omega)}^2$ for some $\rho(\eps)\to 0$. This implies (C1). 

We define now $T_{\epsilon }$ the extension-by-zero operator through the boundary $W$ and $E_{\epsilon}$  the restriction operator to $\Omega$. Note that $T_{\epsilon}$ is well-defined since functions in $W^{2,2}_{0,W}(\Omega)$ vanish on $W$ together with their gradients.  
With these definitions it is straightforward to see that conditions  (C2) and (C3) (i), (ii) are satisfied.

We now prove that condition (C3) (iii) is satisfied. Let $v_{\epsilon }\in   W^{2,2}(\Omega_{\epsilon})\cap W^{1,2}_0(\Omega_{\epsilon})$
be such that $\sup_{\epsilon >0}Q_{\Omega_{\eps}}(v_{\eps})=\sup_{\eps >0}\| v_{\eps}\|^2_{W^{2,2}(\Omega_{\eps})} <\infty $ and let $v\in L^2(\Omega)$ be such that  
$$
{v_{\eps}}_{|_{\Omega}}  \eto v,\ \  {\rm in}\ L^2(\Omega). 
$$
Note  that possibly passing to a subsequence we have that $v_{\eps} \weto v  $ in $W^{2,2}(\Omega )$ and $v_{\eps} \eto v  $ in $W^{1,2}(\Omega )$. Moreover,  proceeding as in the proof of Lemma~\ref{weaklemint} one can show that $v\in W^{1,2}_0(\Omega)$. It remains to prove that $\nabla v=0$ on $W$. Since $v=0$ on $W$ we have that $\frac{\partial v}{\partial x_i}=0$ on $W$ for all $i=1, \dots , N-1$. Thus it suffices to show that $\frac{\partial v}{\partial x_N}=0$ on $W$. For this we  apply Lemma 4.3 from \cite{casado10}.  To do so, for any $i=1, \dots , N-1$ we consider the vector valued function 
\begin{equation}\label{casadovector}V_{\eps}^{(i)}=\left(0, \dots, 0 , -\frac{\partial v_{\eps}}{\partial x_N},0, \dots , 0,   \frac{\partial v_{\eps}}{\partial x_i}\right)
\end{equation}
where the only non-zero entries are the $i$-th and the $N$-th ones. Since $v_{\eps}=0$ on $\partial \Omega_{\eps}$ then $v_{\eps}(\bar x, g_{\eps}(\bar x))=0$ for all $\bar x\in W$. Differentiating this last expression with respect to $ x_i$ for $i=1, \dots , N-1$ we easily get $V_{\eps}^{(i)}\cdot \nu =0 $ on $\Gamma_{\eps}$ where $\nu$ is the normal to the boundary. This allows to apply  Lemma 4.3 from \cite{casado10} to conclude that for any $i=1, \dots , N-1$
$$
\frac{\partial v (\bar x, 0)}{\partial x_N}\frac{\partial b (\bar y)}{\partial y_i}=0,\ \ \   {\rm a.e.\ on}\ W\times Y.
$$
Since $b$ is non-constant, we deduce that $\frac{\partial v}{\partial x_N}=0$ on $W$. Thus $\nabla v=0$ on $W$ hence condition (C) is satisfied. 

The proof of statement (iii) is carried out in the next section. \hfill $\Box$

\section{Critical case $\alpha =3/2$ (proof of Theorem~\ref{main} iii)}
\label{critical}

We divide the proof into several steps organized as subsections. 
There are several important ingredients in the proof of the critical case. The first thing is  to consider again the diffeomorphism $\Phi_\eps:\Omega_\eps\to \Omega$, particularized for this situation and study some of its properties. This is done in  Subsection \ref{SS-diffeomorphism}.  This diffeomorphism will generate its pullback transformation, that we denote again by $T_\eps$, which will allow us to transform functions defined in $\Omega$ to functions defined in $\Omega_\eps$ via composition with the diffeomorphism $\Phi_\eps$, see \eqref{pullback} below. 

With this transformation, we will consider the weak formulation of our problem with a test function of the type $T_\eps\varphi$ with $\varphi$ a test function in $\Omega$, see \eqref{poisson1} below,  and we will be able to easily  pass to the limit in all terms except in a term 
of the type 
$$\int_{W\times(-\eps,0)}D^2v_\eps : D^2{T_\eps \varphi}dx,$$
that requires a deeper analysis.
Notice that this term carries the information of the oscillations of the function $v_\eps$ and the oscillations of the domain, which are coded in the transformation $T_\eps$. Therefore, it is not surprising that this is the the most complicated term to analyze.  We will  treat this term using the unfolding operator method from homogenization. The definition and main ingredients of this tool particularized to our case are contained in Subsection \ref{unfolding}. 

The weak formulation of our problem and passing to the limit in all the terms, including the difficult one is carried out in Subsection \ref{macroscopic}, proving Theorem \ref{mainmacroweak}.  In the limit problem, there is an auxiliary function $\hat v$, which needs to be characterized. As it is customary, this is done by considering another ``oscillatory test function'' and passing to the limit appropriately. The particular calculations to characterize $\hat v$ are contained in Subsection \ref{microscopic}, which concludes the proof of Theorem \ref{main} $iii)$.


\subsection{A special transformation  from $\Omega_{\eps}$ to $\Omega$}
\label{SS-diffeomorphism}

For $\eps>0$ small, we will use the diffeomorphism $\Phi_{\epsilon }: \Omega_{\eps}\to \Omega$ defined  in the proof of Lemma \ref{weaklemint} with $m=2$, $a=-1$,  $g(\cdot)\equiv 0$, and $\tilde g_{\eps}(\cdot)\equiv -\eps$. For the convenience of the reader we write it here explicitly : $\Phi_{\epsilon }(\bar x, x_N)=(\bar x, x_N-h_{\epsilon }(\bar x,x_N))$ for all $(\bar x, x_N)\in \bar\Omega_{\epsilon }$ where
\begin{equation}
\label{accabis}
h_{\epsilon }(\bar x, x_N)=\left\{\begin{array}{ll}0,& \ \ {\rm if }\ -1\le x_N\le -\eps \vspace{2mm}  \\
g_{\epsilon }(\bar x)\left( \frac{x_N+\epsilon }{g_{\epsilon }(\bar x)+\epsilon }   \right)^{3},& \ \ {\rm if }\  -\eps<x_N\le  g_{\epsilon }(\bar x)\, .
\end{array} \right.
\end{equation}

We note  that $\Phi_{\eps}$ is in fact  well-defined for any $\alpha >1$ provided $\eps>0$ is sufficiently small.  For this reason, although this section is devoted to  the case $\alpha =3/2$, we shall try to keep track of $\alpha$ in all formulas and statements where the specific value $\alpha =3/2$ is not required.

The proof of the following lemma follows by straightforward computations. 

\begin{lem} \label{accaest}  The map $\Phi_{\epsilon }$ is a diffeomorphism of class $C^2$ and there exists a constant $c>0$ independent of $\eps$ such that 
$$
|h_{\eps}|\le c \eps^{\alpha},\ \ \left| \frac{\partial h_{\eps}}{\partial x_i}\right|\le c \eps^{\alpha -1},\ \  \left| \frac{\partial^2 h_{\eps}}{\partial x_i\partial x_j}\right|\le c \eps^{\alpha -2},
$$
for all $\eps >0$ sufficiently small. 
\end{lem}

In the sequel we shall use the pull-back operator $T_{\eps}$ associated with $\Phi_{\eps}$. Namely, $T_{\eps}$ is the operator 
\begin{equation}\label{pullback}
\begin{array}{rl}
T_\eps: L^2(\Omega)&\to L^2(\Omega_{\eps}) \\
u&\mapsto u\circ \Phi_{\eps} 
\end{array}
\end{equation}
 Note that $T_{\eps}$ is a linear  homeomorphism and its restrictions to the spaces $W^{1,2}_0(\Omega)$ and $W^{2,2}(\Omega)$ define linear homeomorphisms onto  $W^{1,2}_0(\Omega_{\eps})$ and $W^{2,2}(\Omega_{\eps})$ respectively.  In particular,
$T_{\eps}$ is an isomorphism between the spaces $V(\Omega)=W^{2,2}(\Omega)\cap W^{1,2}_0(\Omega)$ and $V(\Omega_{\eps})=W^{2,2}(\Omega_{\eps})\cap W^{1,2}_0(\Omega_{\eps})$. We also note that for any $\alpha >1$ the operator norm $\|T_{\eps}\|_{\mathcal{L}(W^{1,2}_0(\Omega ), W^{1,2}_0(\Omega_{\eps}))}$ is uniformly bounded with respect to $\eps$, while the operator norm  $\|T_{\eps}\|_{\mathcal{L}(W^{2,2}(\Omega ), W^{2,2}(\Omega_{\eps}))}$ is uniformly bounded with respect to $\eps $  only if $\alpha \geq 2$.

\begin{rem}\label{No-condition-C}
As a difference from Theorem  \ref{main},  case i) and ii) in the critical case $\alpha=3/2$ we will not be able to show condition (C).  We explain here where this condition fails. Notice first, that $T_\eps$ is the natural candidate to show condition (C2), see Definition \ref{conc}.  The other natural operator $E_\eps:V(\Omega_\eps)\to W^{2,2}(\Omega)$ should be the restriction operator.   With respect to $K_\eps$ we have two different and ``natural'' options: $K_\eps=\Omega$ or $K_\eps=W\times (-1,-\eps)$. For both options, condition \eqref{limmeas} and (C1) hold in an easy way.  Moreover, condition (C3) and (C2) iii) also hold. 

The main difficulty is with condition (C2) i) or (C2) ii) depending on the choice of $K_\eps$.  In case $K_\eps=\Omega$, then (C2) i) does not hold and in case $K_\eps=W\times (-1,-\eps)$ then (C2) ii) does not hold. As a matter of fact it will be seen later that $\|T_\eps \varphi\|_{W^{2,2}(W\times (-\eps,0)}$ does not go to 0 for most $\varphi\in V(\Omega)$.     

Hence, it is not possible to show condition (C) for this case. 
\end{rem}

\subsection{Unfolding operator}
\label{unfolding}

We will see that the limiting problem will contain an extra boundary term, which represents the  interplay between  the boundary oscillations and the boundary conditions.   In order to identify the limiting problem in the case $\alpha =3/2$ and prove Theorem~\ref{main} (iii), we shall use the unfolding operator method. In this section we recall the definition of the unfolding operator and some of its properties. We follow the approach of Casado-D\'iaz et al.~\cite{casado, casado10} and we consider an unfolding operator which is an anisotropic version of the classical unfolding operator discussed in Cioranescu, Damlamian and Griso~\cite{ciodamgri}. We note that the well-known  properties of the standard unfolding operator have to be slightly modified. For example, in the  exact integration formula stated in Lemma~\ref{exact} below, an extra factor $\eps$ appears in the right-hand side. Moreover, the limiting function $V$ in Lemma~\ref{unfoldconv} below turns out to be $Y$-periodic  while, keeping in mind in the classical unfolding method,  one would expect that  the  limiting function would be the sum of a periodic function and a polynomial of the second degree in the variables $y$ (cf.,  \cite[Thm.~3.6]{ciodamgri}).

For any $k\in {\mathbb{Z}}^{N-1}$ and $\eps >0$ we consider the $\eps$-cell $C_{\eps}^k=\eps k+\eps Y$, where as above, the basic cell $Y$ is given by $Y=(-\frac{1}{2},\frac{1}{2})^{N-1}$.  Let $I_{W,\eps}=\{k\in {\mathbb{Z}}^{N-1}:\ C_{\eps}^k\subset W\}$. We set 
$$
\widehat W_{\eps}= \bigcup_{k\in I_{W,\eps}} C_{\eps}^k. 
$$

Then we recall the following 

\begin{defn} Let $u$ be a real-valued function defined on $\Omega$. For any $\eps >0$ sufficiently small the unfolding $\hat u$ of $u$ is the real-valued function defined on $\widehat W_{\eps}\times Y\times (-1/\eps ,0)$  by 
\begin{equation}
\hat u (\bar x, \bar y , y_N)= u\left(\eps \left[\frac{\bar x}{\eps} \right]+\eps \bar y, \epsilon y_N   \right),
\end{equation}
for all $(\bar x, \bar y , y_N)\in \widehat W_{\eps}\times Y\times (-1/\eps ,0)$, where $ \left[\frac{\bar x}{\eps} \right]$ denotes the integer part of the vector $\frac{\bar x}{\eps} $ with respect to the unit cell $Y$, that is,$ \left[\frac{\bar x}{\eps} \right]=k\in {\mathbb{Z}}^{N-1}$ if and only if $\bar x\in C_\eps^k$ . 
\end{defn}

We also recall the following  lemma 
the proof of which can be carried out exactly as in the standard case discussed in  \cite{ciodamgri}.

\begin{lem}\label{exact}{\bf (Exact integration formula)} Lat $a\in [-1, 0) $ be fixed. Then  
$$
\int_{\widehat W_{\epsilon}\times (a, 0  ) }u(x)dx=\eps\int_{\widehat W_{\epsilon}\times Y\times (a/ \eps, 0  )  }\hat u (\bar x, y)d\bar xdy
$$
for all $u\in L^1(\Omega)$ and $\eps >0$ sufficiently small.  
\end{lem}

We denote by $W^{2,2}_{Per_Y, loc}(Y\times (-\infty,0))$ the space of functions in $ W^{2,2}_{loc}({\mathbb{R}}^N\times (-\infty , 0)) $ which are also $Y$-periodic in the first  $(N-1)$ variables $\bar y$. As customary, we shall often identify functions in $W^{2,2}_{Per_Y, loc}(Y\times (-\infty,0))$ with their restrictions to $Y\times (-\infty , 0) $.  In the sequel we are going to consider functions in $W^{2,2}_{Per_Y, loc}(Y\times (-\infty,0))$ whose second order weak derivatives are square summable in $Y\times (-\infty ,0)$. For this reason we find it convenient to set 
\begin{eqnarray}\label{energyspace}\lefteqn{
w^{2,2}_{Per_Y}(Y\times (-\infty,0)) =\left\{u  \in W^{2,2}_{Per_Y, loc}(Y\times (-\infty,0)): \right. }\nonumber \\ & & \qquad\qquad\qquad\qquad\qquad\qquad\left.   \ \| D^{\alpha }u\|_{L^2 (Y\times (-\infty,0))}<\infty, \ \forall |\alpha |=2 \right\} \, .
\end{eqnarray}

\begin{lem}\label{unfoldconv}The following statements hold:
\begin{itemize}\item[(i)] Let $v_{\eps}\in W^{2,2}(\Omega) $  with $\| v_{\eps}\|_{W^{2,2}(\Omega)}\le M$  for all $\eps >0$. Let $V_{\epsilon}$ be defined by
\begin{equation}
 V_{\eps}(\bar x, y)=\hat v_{\eps}(\bar x, y)-\int_{Y}\hat v_{\eps}(\bar x, \bar y, 0)d\bar y -\int_{Y}\nabla_y\hat v_{\eps}(\bar x, \bar y, 0)d\bar y\cdot y \end{equation}
 for $(\bar x, y)\in \widehat W_{\epsilon}\times Y\times (-1/ \eps, 0  ) $. 
Then there exists $\hat v\in L^2(W, w^{2,2}_{Per_Y}(Y\times (-\infty,0)))$ such that 
\begin{itemize}
\item[(a)]$\frac{V_{\eps}}{\eps^{3/2}}\weto \hat v$ {\rm and} $\frac{\nabla_yV_{\eps}}{\eps^{3/2}}\weto \nabla_y\hat v$ {\rm in} $L^2(W\times Y\times ({d}, 0))$ {\rm for any} ${d}<0$.
\item[(b)] $\frac{D^{\alpha }_yV_{\eps}}{\eps^{3/2}}=\frac{D^{\alpha }_y\hat v_{\eps}}{\eps^{3/2}}\weto D^{\alpha }_y\hat v$  {\rm in} $L^2(W\times Y\times (-\infty, 0))$, {\rm for any } $|\alpha |=2$,
\end{itemize}
where it is understood that  functions $V_{\eps}, \nabla_y V_{\eps}$ and $D^{\alpha }_yV_{\eps}$ are extended by zero in the whole of $W\times Y\times (-\infty ,0)$  outside their natural domain of definition $\widehat W_{\epsilon}\times Y\times (-1/ \eps, 0  ) $. 
\item[(ii)] Let $v\in W^{1,2}(\Omega)$. Then 
$$
\widehat{(T_{\epsilon }v)_{|\Omega}}\eto v(\bar x , 0),\ \ {\rm in }\ \ L^2(W\times Y\times (-1 ,0))
$$
\end{itemize}
\end{lem}

   {\bf Proof.} We start proving statement (i). It is obvious that $D^{\alpha }_yV_{\eps}=D^{\alpha }_y\hat v_{\eps}$    for any  $|\alpha |=2$.  By Lemma~\ref{exact} and the chain rule it follows that
\begin{eqnarray}
\lefteqn{   \int_{\widehat W_{\epsilon}\times Y\times (-1/ \eps, 0  )  }\left|  \frac{D^{\alpha }_yV_{\eps}}{\eps^{3/2}} 
 \right|^2d\bar xdy  = \int_{\widehat W_{\epsilon}\times Y\times (-1/ \eps, 0  )  }\eps |  \widehat{ D^{\alpha }v_{\eps}} 
 |^2d\bar xdy} \nonumber \\
& &\qquad\qquad\qquad\qquad = \int_{\widehat W_{\epsilon}\times (-1,0)}|D^{\alpha }v_{\eps} |^2dx \le  \int_{\Omega}|D^{\alpha }v_{\eps} |^2dx\le  M^2,
\end{eqnarray}
for all $\eps >0$, hence  $\left\|\frac{D^{\alpha }_yV_{\eps}}{\eps^{3/2}}\right\|_{L^2(W\times Y\times (-\infty , 0))
} $ is uniformly bounded with respect to $\eps$.

Note that the operator defined by $\int_Yv(\bar y ,0)d\bar y+\int_Y\nabla v (\bar y, 0)d\bar y \cdot y$  for functions $v$ in a Sobolev space of  the type $W^{2,2}(Y\times ({d}, 0))$ with ${d}<0$, is a projector on the space of polynomials   of the first degree in $y$.  Thus, we can apply the Poincar\'{e}-Wirtinger inequality and conclude that for any ${d}<0$ there exists $c_{d}>0$ such that
$$
\left\|  \frac{V_{\eps}}{\eps^{3/2}}   \right\|_{L^2(W\times Y\times ({d},0))} , \left\|  \frac{\nabla_yV_{\eps}}{\eps^{3/2}}   \right\|_{L^2(W\times Y\times ({d},0))}\le c_{{d}}\sum_{|\alpha |=2} \left\| \frac{ D^\alpha_yV_{\eps}}{\eps^{3/2}}  \right\|_{L^2(W\times Y\times ({d},0))}\le c_{{d}}M,
$$
for all $\eps >0$. A standard argument implies the existence of a real-valued function $\hat v$ defined on $W\times Y\times (-\infty , 0)$ which admits weak derivatives up to the second order locally in the variable $y$,  such that $\hat v, \nabla_y \hat v \in L^2(W\times Y\times ({d},0))$ for any ${d}<0$, $D^{\alpha }_y\hat v\in L^2(W\times Y\times (-\infty ,0))$,  and   such that statements (a) and (b) hold. 

It remains to prove that $\hat v$ is $Y$-periodic in the variables $\bar y$. Note that 
$$
\int_Y\nabla_y\hat v (\bar x, \bar y, 0) d\bar y =0,
$$
for almost all $x\in W$, hence it suffices to prove that $\nabla_y \hat v $ is $Y$-periodic in the variables $\bar y$.  We note that 
$$
\frac{\nabla_y V_{\eps}(\bar x, y)}{\eps^{3/2}}=\frac{\widehat{\nabla v_{\eps}}(\bar x, y)-\int_Y \widehat{\nabla v_{\eps}}(\bar x, \bar y, 0)d\bar y   }{\sqrt{\eps}}.
$$
Thus, in order conclude it is simply enough to apply the same argument in Step 3 of the proof of Lemma~4.3 in Casado-D\'iaz et al.\cite{casado10}  to the function $\nabla v_{\eps}$.

We now prove statement (ii). If  $v\in C^{\infty }(\bar \Omega)$, we easily can see that 
\begin{eqnarray}\lefteqn{
\int_{\widehat{W_{\eps}}\times Y\times (-1,0)}\left| \widehat{(T_{\epsilon }v)_{|\Omega}}- v(\bar x , 0) \right|^2d\bar x dy=
\int_{Y\times (-1,0)}\int_{\cup_{k\in I_{W,\eps}}C_{\eps}^k }\left| \widehat{T_{\epsilon }v}- v(\bar x , 0) \right|^2d\bar x dy
}\nonumber  \\ & & 
 = \int_{-1}^0\sum_{k\in I_{W,\eps}}\int_{C_{\eps}^k}\int_Y\left|  (T_{\eps}v)(\eps\left[\frac{\bar x}{\eps} \right]+\eps \bar y, \eps y_N   )  - v(\bar x , 0) \right|^2d\bar yd\bar x dy_N \nonumber
\\ & & 
= \int_{-1}^0\sum_{k\in I_{W,\eps}}\int_{C_{\eps}^k}\int_{C_{\eps}^k}\left|  (T_{\eps}v)( \bar z, \eps y_N   )  - v(\bar x , 0) \right|^2\frac{d\bar z}{\eps^{N-1}}d\bar x dy_N \nonumber
\\ & & 
= \int_{-1}^0\sum_{k\in I_{W,\eps}}\int_{C_{\eps}^k}\int_{C_{\eps}^k}\left|  v( \bar z, \eps y_N -h_{\eps}(\bar z, \eps y_N)  )  - v(\bar x , 0) \right|^2\frac{d\bar z}{\eps^{N-1}}d\bar x dy_N \nonumber
\\ & & 
\le c\int_{-1}^0\sum_{k\in I_{W,\eps}}\int_{C_{\eps}^k}\int_{C_{\eps}^k}\left| \bar z -\bar x \right|^2+ \left| \eps y_N-h_{\eps}(\bar z, \eps y_N) \right|^2\frac{d\bar z}{\eps^{N-1}}d\bar x dy_N \nonumber
\\ & & 
\le c\int_{-1}^0\sum_{k\in I_{W,\eps}}\int_{C_{\eps}^k}\eps^2d\bar x dy_N\le c \eps^2,
\end{eqnarray}
hence statement (ii) is proved for smooth functions.  In the case of an arbitrary function $v\in W^{1,2}(\Omega)$, we use an approximation argument. Namely, we consider a sequence 
$v_n\in C^{\infty }(\bar \Omega )$ converging to $v$ in $W^{1,2}(\Omega)$ as $n\to \infty$ and we note that 
\begin{equation}\label{tria}\|\widehat{(T_{\epsilon }v)_{|\Omega}}-v(x,0)\|\le 
\|
\widehat{(T_{\epsilon }v)_{|\Omega}}-\widehat{(T_{\epsilon }v_n)_{|\Omega}}\|+\|\widehat{(T_{\epsilon }v_n)_{|\Omega}}-v_n(\bar x,0)\|+\|v_n(\bar x, 0)-v(\bar x, 0)\|
\end{equation}
where all norms are taken in $ L^2(W\times Y\times (-1 ,0))$. Since statement (ii) holds for smooth functions, the second term in the right hand-side of (\ref{tria}) goes to zero as $\eps \to 0$. Moreover, by the continuity of the trace operator, also the third term in the right hand-side of (\ref{tria}) goes to zero as $\eps \to 0$. We now consider the first term  in the right hand-side of (\ref{tria}) . By Lemma \ref{exact} and changing variables in integrals, we have
\begin{eqnarray}\label{tria1}\lefteqn{
\int_{\widehat{W_{\eps}}\times Y\times (-1,0)}\left| \widehat{T_{\epsilon }v}- \widehat{T_{\epsilon }v_n} \right|^2d\bar x dy}\nonumber \\ & &  =
\eps^{-1}\int_{\widehat{W_{\eps}}\times (-\eps,0)}|v\circ \Phi_{\eps}-v_n\circ \Phi_{\eps}|^2dx\nonumber \\ & & 
\le c 
\eps^{-1}\int_{\Phi_{\eps}(\widehat{W_{\eps}}\times (-\eps,0))}|v-v_n|^2dx\le c\| v-v_n\|^2_{W^{1,2}(\Omega )}
\end{eqnarray}
where the last inequality is deduced by the fact that the diameter in the $x_N$-direction of the set $\Phi_{\eps}(\widehat{W_{\eps}}\times (-\eps,0))$
is $O(\eps)$ as $\eps\to 0$ and that the function $|v-v_n|$ is bounded in almost all vertical lines. 
Since the right hand-side of (\ref{tria1}) goes to zero as $n\to \infty$, we easily conclude. \hfill $\Box$

\subsection{Weak macroscopic limiting problem}
\label{macroscopic}

 Let $f_{\eps}\in L^{2}(\Omega_{\eps}) $ and $f\in  L^{2}(\Omega)$ be such that $f_{\eps}\weto f$  in $L^2({\mathbb{R}}^N)$ with the understanding that such functions are extended by zero outside $\Omega_{\eps}$ and $\Omega $ respectively.  Let $v_{\eps}\in V(\Omega_{\eps})=W^{2,2}(\Omega_{\eps})\cap W^{1,2}_0(\Omega_{\eps})$  be such that 
\begin{equation}\label{poisson}
H_{\Omega_{\eps}, I}v_{\eps}=f_{\eps}
\end{equation}
for all $\eps > 0$ small enough. By (\ref{poisson}) it follows that $\|v_{\eps}\|_{ W^{2,2}(\Omega_{\eps})}\le M$ for all $\eps >0$ sufficiently small hence, possibly passing to a subsequence, there exists $v\in W^{2,2}(\Omega )\cap W^{1,2}_0(\Omega )$ such that $v_{\eps}\weto v$ in $ W^{2,2}(\Omega )$ and  $v_{\eps}\eto v$ in $ L^{2}(\Omega)$. 

Let $\varphi \in V(\Omega )= W^{2,2}(\Omega )\cap W^{1,2}_0(\Omega ) $ be a  fixed test function. Since $T_{\eps}\varphi \in V(\Omega_{\eps})$, by (\ref{poisson}) it follows that
\begin{equation}
\label{poisson1}
\int_{\Omega_{\eps}}D^2v_{\eps}: D^2T_{\eps}\varphi dx +\int_{\Omega_{\eps}}v_{\eps}T_{\eps}\varphi dx =\int_{\Omega_{\eps}}f_{\eps}T_{\eps}\varphi dx.
\end{equation} 

It is easy to see that 
\begin{equation}\label{poisson2}
\int_{\Omega_{\eps}}v_{\eps}T_{\eps}\varphi dx\eto \int_{\Omega}v\varphi dx,\ \ {\rm and}\ \ \int_{\Omega_{\eps}}f_{\eps}T_{\eps}\varphi dx\eto \int_{\Omega}f\varphi dx.
\end{equation}

We now consider the first term in the left hand-side of (\ref{poisson1}).
It is convenient to set $K_{\eps}= W\times (-1, -\eps)$  so that 
\begin{equation}\label{poisson3}
\int_{\Omega_{\eps}}D^2v_{\eps}: D^2T_{\eps}\varphi dx = \int_{\Omega_{\eps}\setminus \Omega }D^2v_{\eps}: D^2T_{\eps}\varphi dx+ \int_{\Omega\setminus K_{\eps}}D^2v_{\eps}: D^2T_{\eps}\varphi dx + \int_{K_{\eps}}D^2v_{\eps}: D^2T_{\eps}\varphi dx
\end{equation}
Since $D^2T_{\eps}\varphi=D^2\varphi $ in $K_{\eps}$ we have that 
\begin{equation}\label{poisson4}
\int_{K_{\eps}}D^2v_{\eps}: D^2T_{\eps}\varphi dx=\int_{K_{\eps}}D^2v_{\eps}: D^2\varphi dx \eto \int_{\Omega}D^2v: D^2\varphi dx.
\end{equation}


Moreover, one can prove that
\begin{equation}
\label{poisson5}
 \int_{\Omega_{\eps}\setminus \Omega }D^2v_{\eps}: D^2T_{\eps}\varphi dx\eto 0 .
\end{equation}
Indeed, by changing  variables in integrals, using  the chain rule and  Lemma~\ref{accaest} we get the following inequalities (here and in the sequel, to shorten notation we drop the summation symbols): 
\begin{eqnarray}\label{poisson51}\lefteqn{\left(
 \int_{\Omega_{\eps}\setminus \Omega }D^2v_{\eps}: D^2T_{\eps}\varphi dx\right)^2\le c   \int_{\Omega_{\eps}\setminus \Omega }\left| \frac{\partial^2 T_{\eps}\varphi}{\partial x_i\partial x_j}   \right|^2dx}\nonumber \\
& & 
\le c   \int_{\Omega_{\eps}\setminus \Omega }\biggl| \frac{\partial^2 \varphi  }{\partial x_k\partial x_l} (\Phi_{\epsilon}(x))\frac{\partial \Phi_{\eps}^{(k)}}{\partial x_i} \frac{\partial \Phi_{\eps}^{(l)}}{\partial x_j}  \biggr|^2 +  \biggl|\frac{\partial \varphi }{\partial x_k}(\Phi_{\epsilon}(x))   \frac{\partial^2   (\Phi^{(k)}_{\epsilon}(x)) }{\partial x_i\partial x_j} \biggr|^2  dx  \nonumber \\
& & \le  c  \int_{\Omega_{\eps}\setminus \Omega }\biggl| \frac{\partial^2 \varphi  }{\partial x_k\partial x_l} (\Phi_{\epsilon}(x))\biggr|^2dx +  c \eps^{-1}  \int_{\Omega_{\eps}\setminus \Omega } \biggl|\frac{\partial \varphi }{\partial x_k}(\Phi_{\epsilon}(x)) \biggr|^2  dx \nonumber\\
& & \le  c  \int_{\Phi_{\eps}(\Omega_{\eps}\setminus \Omega )}\biggl| \frac{\partial^2 \varphi  }{\partial x_k\partial x_l} \biggr|^2dx +  c \eps^{-1}  \int_{\Phi_{\eps}(\Omega_{\eps}\setminus \Omega )} \biggl|\frac{\partial \varphi }{\partial x_k} \biggr|^2  dx  .
\end{eqnarray}
Observe that 
\begin{eqnarray}
\Phi_\eps(\Omega_\eps\setminus\Omega)&=&\{(\bar x,x_N-h_\eps(\bar x,x_N)): \bar x\in W,\ 0<x_N<g_\eps(\bar x)\} \nonumber \\
&\subset& \{(\bar x,z_N): \bar x\in W, -g_\eps(\bar x)<z_N<0\}\nonumber \\
&\subset& \{(\bar x,z_N): \bar x\in W, -\eps^{3/2}b_0<z_N<0\},
\end{eqnarray}
where $b_0=\|b(\cdot)\|_{L^\infty(W)}$.  Hence, $|\Phi_\eps(\Omega_\eps\setminus\Omega)|\leq c\eps^{3/2}\eto 0$ and therefore
$$ \int_{\Phi_{\eps}(\Omega_{\eps}\setminus \Omega )}\biggl| \frac{\partial^2 \varphi  }{\partial x_k\partial x_l} \biggr|^2dx\eto 0.$$

Moreover,  notice that $\frac{\partial \varphi }{\partial x_k}(\bar x, \cdot)\in W^{1,2}(-1,0)$, a.e. $\bar x\in W$ and therefore 
$$\left\|\frac{\partial \varphi }{\partial x_k}(\bar x, \cdot)\right\|_{L^\infty(-1,0)}^2\leq C\left\|\frac{\partial \varphi }{\partial x_k}(\bar x, \cdot)\right\|_{W^{1,2}(-1,0)}^2 \hbox{ a.e. }\bar x\in W.$$

Hence, the last term from \eqref{poisson51} is analyzed as follows
\begin{eqnarray}\lefteqn{ c\eps^{-1}\int_{\Phi_{\eps}(\Omega_{\eps}\setminus \Omega )} \biggl|\frac{\partial \varphi }{\partial x_k} \biggr|^2  dx\leq c\eps^{-1}\eps^{3/2}b_0\int_{W}\bigg\|\frac{\partial \varphi }{\partial x_k}(\bar x, \cdot)\bigg\|_{L^\infty(-1,0)}^2 d\bar x }     \nonumber \\
& & \qquad\leq C\eps^{1/2}\int_{W}\bigg\|\frac{\partial \varphi }{\partial x_k}(\bar x, \cdot )\bigg\|^2_{W^{1,2}(-1,0)}d\bar x\leq C\eps^{1/2}\|\varphi\|_{W^{2,2}(\Omega)}^2\eto 0.
\end{eqnarray}

%

We now consider the second term in the right hand-side of (\ref{poisson3}). It is convenient to set $Q_{\eps}=\widehat W_{\eps}\times (-\eps , 0)$ so that 
\begin{equation}
\label{poisson6}
\int_{\Omega\setminus K_{\eps}}D^2v_{\eps}: D^2T_{\eps}\varphi dx = \int_{Q_{\eps} }D^2v_{\eps}: D^2T_{\eps}\varphi dx +\int_{\Omega\setminus (K_{\eps}\cup Q_{\eps}) }D^2v_{\eps}: D^2T_{\eps}\varphi dx .
\end{equation}


One can prove that
\begin{equation}
\label{poisson7}
\int_{\Omega\setminus (K_{\eps}\cup Q_{\eps}) }D^2v_{\eps}: D^2T_{\eps}\varphi dx\eto 0.
\end{equation}
Indeed, proceeding exactly as in (\ref{poisson51}) we get 
\begin{eqnarray}\label{poisson52}\lefteqn{\left(
 \int_{\Omega\setminus (K_{\eps}\cup Q_{\eps}) }D^2v_{\eps}: D^2T_{\eps}\varphi dx\right)^2}\nonumber \\
& &  \le  c  \int_{\Phi_{\eps}(  \Omega\setminus (K_{\eps}\cup Q_{\eps})}\biggl| \frac{\partial^2 \varphi  }{\partial x_k\partial x_l} \biggr|^2dx +  c \eps^{-1}  \int_{\Phi_{\eps}(\Omega\setminus (K_{\eps}\cup Q_{\eps}) )} \biggl|\frac{\partial \varphi }{\partial x_k} \biggr|^2  dx\nonumber \\
& &  \le  c  \int_{\Phi_{\eps}( \Omega\setminus (K_{\eps}\cup Q_{\eps} )}\biggl| \frac{\partial^2 \varphi  }{\partial x_k\partial x_l} \biggr|^2dx +  c \| \varphi   \|^2_{W^{2,2}(\Omega \setminus (\widehat{W_{\eps}}\times (-1,0)   )   )} ,
\end{eqnarray}
where now we have used the fact that the diameter of the set $\Phi_{\eps}( \Omega\setminus (K_{\eps}\cup Q_{\eps} )$ in the direction $x_N$ is $O(\eps)$ as $\eps \to 0$.  Moreover, since $|\Omega\setminus (\widehat{W_{\eps}}\times (-1,0))|\to 0$ and $\varphi$ is a fixed function,  we have 
$\| \varphi   \|^2_{W^{2,2}(\Omega \setminus (\widehat{W_{\eps}}\times (-1,0)   )   )} \eto 0$. Thus, (\ref{poisson7}) follows.

\par\medskip

It remains to analyze the first term  in the right hand-side of (\ref{poisson6}). 
To do so we need the following technical lemma.

\begin{lem}\label{easy} For all $y\in Y\times (-1,0)$ and $i,j=1,\dots , N$, the functions $\hat h_{\eps} (\bar x , y)$, $\widehat{ \frac{\partial h_{\eps} }{\partial x_i }  } (\bar x , y)$ and $\widehat{ \frac{\partial^2h_{\eps} }{\partial x_i\partial x_j }  } (\bar x , y) $ are independent of $\bar x$. Moreover, $\hat h_{\eps} (\bar x , y)=O(\eps^{3/2})$, $\widehat{ \frac{\partial h_{\eps} }{\partial x_i }  } (\bar x , y)=O(\eps^{1/2})$ as $\eps\to 0$  and 
\begin{equation}\label{easy0} \epsilon^{\frac{1}{2}}
\widehat{ \frac{\partial^2h_{\eps} }{\partial x_i\partial x_j }  } (\bar x , y)  \eto \frac{\partial^2 \left( b(\bar y)  (y_N+1)^3\right)
  }{\partial y_i\partial y_j},
\end{equation}
for all $i,j=1, \dots , N$, 
uniformly in $y\in Y\times (-1,0) $.
\end{lem}

{\bf Proof.} The independence of the functions in the statement   from  $\bar x$ is easily deduced by the periodicity of the function 
$b$ and  the definition of $h_\eps$ in (\ref{accabis}). The rest of the proof follows by straightforward computations and   we report only those required for the proof of (\ref{easy0}) for the convenience of the reader. Note that 
\begin{eqnarray}\label{easy1}
\widehat{ \frac{\partial^2h_{\eps} }{\partial x_i\partial x_j }  }& = &\eps^{\alpha-2}\frac{\partial^2 b }{\partial y_i\partial y_j }\widehat{\left( \frac{x_N+\eps}{g_{\eps}(\bar x)+\eps}  \right)^3  } +\eps^{\alpha -1}\frac{\partial b}{\partial y_i}\widehat{\frac{\partial}{\partial x_j}}\left( \frac{x_N+\eps}{g_{\eps}(\bar x)+\eps}  \right)^3\nonumber \\
& & +\eps^{\alpha -1}\frac{\partial b}{\partial y_j}\widehat{\frac{\partial}{\partial x_i}}\left( \frac{x_N+\eps}{g_{\eps}(\bar x)+\eps}  \right)^3+\eps^{\alpha}b(\bar y)\widehat{\frac{\partial^2}{\partial x_i\partial x_j}}\left( \frac{x_N+\eps}{g_{\eps}(\bar x)+\eps}  \right)^3 .
\end{eqnarray}
Moreover, 
\begin{equation}\label{easy2}
\widehat{\frac{\partial}{\partial x_i}}\left( \frac{x_N+\eps}{g_{\eps}(\bar x)+\eps}  \right)^3 = \frac{ 3\delta_{iN}\eps^{-1}(y_N+1)^2   }{   ( \epsilon^{\alpha -1}b(\bar y)+1   )^3  }
-\frac{3\delta_{iN}\epsilon^{\alpha-2}(y_N+1)^3\frac{\partial b}{\partial y_i}}{(\eps^{\alpha -1}b(\bar y)+1  )^4}
\end{equation}
and
\begin{eqnarray}\label{easy3}\lefteqn{
\widehat{\frac{\partial^2}{\partial x_i\partial x_j}}\left( \frac{x_N+\eps}{g_{\eps}(\bar x)+\eps}  \right)^3 = \frac{6\delta_{iN}\delta_{jN}\eps^{-2}(y_N+1)   }{ (\eps^{\alpha -1}b(\bar y)+1  )^3 }  } \nonumber \\
& & -\frac{9 \eps^{\alpha -3}(y_N+1)^2 (\delta_{iN}\frac{\partial b}{\partial y_j}+ \delta_{jN}\frac{\partial b}{\partial y_i} )  }{ (\eps^{\alpha -1}b(\bar y)+1  )^4 }+\frac{12 \eps^{2\alpha -4}(y_N+1)^3 \frac{\partial b}{\partial y_i}\frac{\partial b}{\partial y_j}   }{ (\eps^{\alpha -1}b(\bar y)+1  )^5 }\nonumber  \\
& & -\frac{3 \eps^{\alpha -3}(y_N+1)^3 \frac{\partial^2 b}{\partial y_i\partial y_i}   }{ (\eps^{\alpha -1}b(\bar y)+1  )^4 }
\end{eqnarray}

By combining (\ref{easy1})-(\ref{easy3}) we easily get (\ref{easy0}). \hfill $\Box$\\


We are now ready to prove the following


\begin{lem}\label{macrolem}
Let  $\hat v\in  L^2(W, w^{2,2}_{per_y}(Y\times (-\infty , 0)))$ be as in Lemma~\ref{unfoldconv}. Then 
\begin{equation}\label{macrolem0}
 \int_{Q_{\eps} }D^2v_{\eps}: D^2T_{\eps}\varphi dx \eto -  \int_{W} \int_{ Y\times (-1,0)} D^2_y \hat v(\bar x, y) : D^2_{y}( b(\bar y)(1+y_N)^3  ) dy   \frac{\partial \varphi}{\partial x_N}(\bar x, 0)d\bar x,
\end{equation}
where $D^2_{y}$ denotes the Hessian matrix in the variable $y$. 
\end{lem}

 {\bf Proof.}
By Lemma~\ref{exact} and the chain rule we get
\begin{eqnarray}\label{macrolem1} \lefteqn{
 \int_{Q_{\eps} }D^2v_{\eps}: D^2T_{\eps}\varphi dx}\nonumber \\ 
& & = \eps \int_{\widehat{W_{\eps}}\times Y\times (-1,0)}\widehat{D^2v_{\eps}}: \widehat{D^2T_{\eps}\varphi}d\bar xdy
=\eps^{-3} \int_{\widehat{W_{\eps}}\times Y\times (-1,0)}  D^2_y\widehat{v_{\eps}}: D^2_y\widehat{T_{\eps}\varphi}   d\bar xdy\nonumber \\ 
& & = \eps^{-3}\int_{\widehat{W_{\eps}}\times Y\times (-1,0) }\frac{\partial^2 \widehat{v_{\epsilon}} }{\partial y_i\partial y_j}
\frac{\partial^2 \varphi  }{\partial x_k\partial x_l} (\widehat\Phi_{\epsilon}(y))\frac{\partial \widehat\Phi_{\eps}^{(k)}}{\partial y_i} \frac{\partial\widehat \Phi_{\eps}^{(l)}}{\partial y_j} d\bar xdy \label{macrolem1a} \\
& & + \eps^{-3}\int_{\widehat{W_{\eps}}\times Y\times (-1,0) }\frac{\partial^2 \widehat{v_{\epsilon}} }{\partial y_i\partial y_j}
\frac{\partial \varphi  }{\partial x_k} (\widehat\Phi_{\epsilon}(y)) \frac{\partial^2 \widehat\Phi_{\eps}^{(k)}}{\partial y_i\partial y_j}  d\bar xdy . \label{macrolem1b}\end{eqnarray} 

We note that 
\begin{equation}\label{macrolem2} 
\frac{\partial\widehat\Phi_{\eps}^{(k)}}{\partial y_i}=\left\{ 
\begin{array}{ll}\eps \delta_{ki},\ \ & {\rm if}\ k\ne N,\\
\eps\delta_{Ni}-\eps\widehat{ \frac{\partial h_{\eps}}{\partial x_i} },\ \ & {\rm if}\ k= N.
\end{array}
\right.
\end{equation}

By Lemma~\ref{unfoldconv}  (i), we have that $\| \frac{\partial^2 \widehat{v_{\epsilon}} }{\partial y_i\partial y_j}\|_{L^2(\widehat{W_{\eps}}\times Y\times (-1,0)) }=O(\eps^{3/2})$ as $\eps \to 0$
and by (\ref{macrolem2}) we have that $\frac{\partial\widehat\Phi_{\eps}^{(k)}}{\partial y_i}=O(\eps)$  as $\eps \to 0$. Using also Lemma~\ref{exact} once more, we get
\begin{eqnarray}\label{macrolem3}\lefteqn{\left|
\eps^{-3}\int_{\widehat{W_{\eps}}   \times Y\times (-1,0) }\frac{\partial^2 \widehat{v_{\epsilon}} }{\partial y_i\partial y_j}
\frac{\partial^2 \varphi  }{\partial x_k\partial x_l} (\widehat\Phi_{\epsilon}(y))\frac{\partial \widehat\Phi_{\eps}^{(k)}}{\partial y_i} \frac{\partial\widehat \Phi_{\eps}^{(l)}}{\partial y_j} d\bar xdy\right|}\nonumber \\  & & 
\le c\epsilon^{2} \left\| \eps^{-3/2} \frac{\partial^2 \widehat{v_{\epsilon}} }{\partial y_i\partial y_j}\right\|_{L^2(\widehat{W_{\eps}}\times Y\times (-1,0) )}
\left\|   \eps^{-3/2}  \frac{\partial^2 \varphi  }{\partial x_k\partial x_l} (\widehat\Phi_{\epsilon}(y))\right\|_{L^2(\widehat{W_{\eps}}\times Y\times (-1,0) )}\nonumber \\
& & 
\le c \eps^{1/2}\left\|   \frac{\partial^2 \varphi  }{\partial x_k\partial x_l} (\widehat\Phi_{\epsilon}(y))\right\|_{L^2(\widehat{W_{\eps}}\times Y\times (-1,0) )}
\le  c   \left\|  \frac{\partial^2 \varphi  }{\partial x_k\partial x_l} (\Phi_{\epsilon}(x))\right\|_{L^2(Q_{\eps})}\nonumber \\
& & 
\le 
c\left\|  \frac{\partial^2 \varphi  }{\partial x_k\partial x_l}  \right\|_{ L^2(\Phi_{\eps}(Q_{\eps}))}.
\end{eqnarray}
Hence,  the integral in (\ref{macrolem1a}) vanishes as $\eps \to 0$.

We note that 
\begin{equation}\label{macrolem4} 
\frac{\partial^2\widehat\Phi_{\eps}^{(k)}}{\partial y_i\partial y_j}=\left\{ 
\begin{array}{ll}0,\ \ & {\rm if}\ k\ne N,\\
-\eps^2\widehat{ \frac{\partial^2 h_{\eps}}{\partial x_i\partial x_j} },\ \ & {\rm if}\ k= N.
\end{array}
\right.
\end{equation}
Therefore, in the sum of integrals in (\ref{macrolem1b})  the terms corresponding to the indexes $k\ne N$ vanish. It remains to analyze the term in (\ref{macrolem1b}) with $k=N$.  

By (\ref{macrolem4}) we  rewrite such term in the form
\begin{eqnarray}\label{macrolem5}\lefteqn{
\eps^{-3}\int_{\widehat{W_{\eps}}\times Y\times (-1,0) }\frac{\partial^2 \widehat{v_{\epsilon}} }{\partial y_i\partial y_j}
\frac{\partial \varphi  }{\partial x_N} (\widehat\Phi_{\epsilon}(y)) \frac{\partial^2 \widehat\Phi_{\eps}^{(N)}}{\partial y_i\partial y_j}  d\bar xdy }\nonumber \\ & & \ \ \ \ 
=-\int_{\widehat{W_{\eps}}\times Y\times (-1,0) }\left( \eps^{-3/2}\frac{\partial^2 \widehat{v_{\epsilon}} }{\partial y_i\partial y_j}\right)
\widehat{T_{\eps}\frac{\partial \varphi}{\partial x_N}   }
  \left(  \eps^{1/2} 
\widehat{  \frac{\partial^2 h_{\eps}}{\partial x_i\partial x_j}}     
\right) d\bar xdy  .
\end{eqnarray}
By applying Lemma~\ref{unfoldconv} (i) to the sequence $v_{\eps}$,   Lemma~\ref{unfoldconv} (ii) with $v$ replaced by $\frac{\partial \varphi}{\partial x_N}$, and  using (\ref{easy0}) we easily deduce that the integral in the right hand-side of (\ref{macrolem5}) converges to the integral in the right-hand side of (\ref{macrolem0}) as $\eps \to 0$. \hfill $\Box$\\


Thus we have proved the following

\begin{thm}\label{mainmacroweak} Let $f_{\eps}\in L^2(\Omega_{\eps})$ and $f\in L^2(\Omega)$ be such that $f_{\eps}\weto f$ in $L^2(\Omega)$. Let $v_{\eps}\in W^{2,2}(\Omega_{\eps})\cap W^{1,2}_0(\Omega_{\eps})$ be the solutions to (\ref{poisson}). Then possibly passing to a subsequence, there exists $v\in  W^{2,2}(\Omega)\cap W^{1,2}_0(\Omega)$ and $\hat v \in L^2(W, w^{2,2}_{Per_Y}(Y\times (-\infty , 0)))$ such that 
$v_{\eps}\weto v$ in $W^{2,2}(\Omega)$, $v_{\eps}\eto v$ in $L^{2}(\Omega)$ and such that statements (a) and (b) in Lemma~\ref{unfoldconv} hold, and such that 
\begin{eqnarray}\label{mainmacroweak1}\lefteqn{
\int_{\Omega}D^2v:D^2\varphi +u\varphi dx  -  \int_{W} \int_{ Y\times (-1,0)} D^2_y \hat v(\bar x, y) : D^2_{y}( b(\bar y)(1+y_N)^3  ) dy   \frac{\partial \varphi}{\partial x_N}(\bar x, 0)d\bar x }\nonumber \\ & & \qquad\qquad\qquad\qquad \qquad\qquad\qquad\qquad \qquad\qquad\qquad\qquad= \int_{\Omega}f\varphi dx ,
\end{eqnarray}
for all $\varphi \in W^{2,2}(\Omega)\cap W^{1,2}_0(\Omega)$.
\end{thm}

\subsection{Characterization of $\hat v$ via a weak microscopic problem}
\label{microscopic}
In this section we plan to characterize the function $\hat v$ defined in Theorem~\ref{mainmacroweak}.  

Let $\psi \in C^{\infty}(\overline{W}\times\overline{Y}\times ]-\infty , 0]    )$
 be such that 
${\rm supp}\, \psi\subset C\times \overline{Y}\times [d,0]$ for some compact set $C\subset W$ and $d\in ]-\infty ,0[$ and such that $\psi (\bar x, \bar y, 0)=0$ for  all $(\bar x, \bar y)\in W\times Y$.  Assume also that $\psi$ is $Y$-periodic in the second $(N-1)$-variables $\bar y$.   We set

\begin{equation}
\label{micro0}
\psi_{\eps}(x)=\eps^{\frac{3}{2}}\psi \left(\bar x, \frac{\bar x}{\eps}, \frac{ x_N}{\eps}\right), 
\end{equation}
for all $\epsilon >0$, $x\in W\times ]-\infty , 0] $. Note that for $\eps $ sufficiently small we have that ${\rm supp }\, \psi_{\eps}\subset \Omega $ and  $\psi_{\eps}\in V(\Omega)$, hence $T_{\eps}\psi_{\eps}$ belongs to $ V(\Omega_{\eps})$ and can be used as test function in the weak formulation of our problem in $\Omega_{\eps}$, that is,   

\begin{equation}
\label{micro1}
\int_{\Omega_{\eps}}D^2v_{\eps}: D^2T_{\eps}\psi_{\eps} dx +\int_{\Omega_{\eps}}v_{\eps}T_{\eps}\psi_{\eps} dx =\int_{\Omega_{\eps}}f_{\eps}T_{\eps}\psi_{\eps} dx\, .
\end{equation} 

By the presence of the factor $\eps^{\frac{3}{2}}$ in (\ref{micro0}), it is easy to see that  
\begin{equation}
\label{micro2}
\int_{\Omega_{\eps}}v_{\eps}T_{\eps}\psi_{\eps} dx\eto 0,\ \ {\rm and}\ \ \int_{\Omega_{\eps}}f_{\eps}T_{\eps}\psi_{\eps} dx\eto 0. 
\end{equation} 
We now consider the first term in the left hand-side of (\ref{micro1}) and we write it in the form
\begin{equation}
\label{micro3}
\int_{\Omega_{\eps}}D^2v_{\eps}: D^2T_{\eps}\psi_{\eps} dx= \int_{\Omega_{\eps}\setminus \Omega}D^2v_{\eps} : D^2T_{\eps}\psi_{\eps} dx+ \int_{\Omega}D^2v_{\eps} :  D^2T_{\eps}\psi_{\eps} dx
\end{equation}


With very similar arguments as the ones used to show (\ref{poisson5}) we can also prove that

\begin{equation}\label{micro3,5}
 \int_{\Omega_{\eps}\setminus \Omega }D^2v_{\eps}: D^2T_{\eps}\psi_\eps dx\eto 0. 
\end{equation}

We can show now the following

\begin{lem}\label{Lemma-micro2}
If $\psi_{\eps}$ is as in (\ref{micro0}) and $v_{\eps}$,  $\hat v$ are  the functions from Theorem~\ref{mainmacroweak}, then
\begin{equation}\label{micro4}
\int_{\Omega}D^2v_{\eps}: D^2T_{\eps}\psi_\eps dx\to \int_{W\times Y\times (-\infty,0)}D^2_y\hat v(\bar x, y) : D^2_y\psi(\bar x, y) d\bar xdy.
\end{equation}

\end{lem}
 {\bf Proof.}    
 First of all we note that by the periodicity of $\psi $ we have that  
\begin{equation}\label{two1}
\widehat {T_{\eps}\psi_{\eps}} (\bar x, y)=\eps^{3/2}\psi \left( \eps \left[\frac{ \bar x}{\eps }  \right]+\eps \bar y , \bar y , y_N-\eps^{-1}h_{\eps}\left( \eps \left[\frac{ \bar x}{\eps }  \right]+\eps \bar y, \eps y_N \right) \right)\, .
\end{equation}
We also note explicitly that 
\begin{equation}
h_{\eps}\left( \eps \left[\frac{ \bar x}{\eps }  \right]+\eps \bar y, \eps y_N \right)=
\left\{\begin{array}{ll}\frac{\eps^{3/2}b(\bar y)(y_N+1)^3}{(\eps^{1/2}b(\bar y)+1 )^3},\ & {\rm if}\ -1\le y_N<0,\nonumber\\
0,\ & {\rm if}\ -1/\eps <y_N<-1,  
\end{array}\right.
\end{equation}
hence 
\begin{equation}\label{two2}
D^{\alpha }_y\left(\eps^{-1}h_{\eps}\left( \eps \left[\frac{ \bar x}{\eps }  \right]+\eps \bar y, \eps y_N \right)\right)\eto 0, \ \ {\rm uniformly\ on}\ W\times Y\times ]-\infty , 0],
\end{equation}
for all $|\alpha |\le 2$.

Since $\psi$ is smooth and has compact support, it is Lipschitz continuous together with its derivatives and it easily follows that  
\begin{equation}\label{two3}\left\|
(D^{\beta }\psi ) \left( \eps \left[\frac{ \bar x}{\eps }  \right]+\eps \bar y , \bar y , y_N-\eps^{-1}h_{\eps}\left( \eps \left[\frac{ \bar x}{\eps }  \right]+\eps \bar y, \eps y_N \right) \right)- D^{\beta }\psi (\bar x, y)\right\|_{L^2(\widehat W_{\eps}\times Y\times ]-\infty , 0[)}\eto 0   ,
\end{equation}
 for any $|\beta |\le 2$. In fact, the square of the norm in (\ref{two3}) can be estimated by
$$
\int_{\widehat W_{\eps}\times Y\times ]-\tilde d , 0[) }\left|  \eps \left[\frac{ \bar x}{\eps }  \right]+\eps \bar y-\bar x  \right|^2+\left|     \eps^{-1}h_{\eps}\left( \eps \left[\frac{ \bar x}{\eps }  \right]+\eps \bar y, \eps y_N \right) \right|^2d\bar xdy
$$
which is clearly $O(\eps)$ as $\eps \to 0$.

By combining (\ref{two1})-(\ref{two3}) and using the chain rule, we get that 
\begin{equation}\label{two4}\eps^{-3/2}
D^{\gamma}_y\widehat {T_{\eps}\psi_{\eps}} (\bar x, y)\eto D_y^{\gamma}\psi (\bar x ,y)
\end{equation}
in $L^2( W\times Y\times ]-\infty , 0[)$, for all $|\gamma |=2$. 

By Lemmas~\ref{exact}, \ref{unfoldconv} and (\ref{two4}),  we conclude that 
\begin{eqnarray}\lefteqn{
\int_{\widehat{W}_{\eps}\times ]-1,0[   }D^2v_{\eps}: D^2T_{\eps}\psi_{\eps}dx }\nonumber \\
& & =\int_{\widehat{W}_{\eps}\times Y\times ]-1/\eps , 0[}\frac{D^2_y\hat v_{\eps}}{\eps^{3/2}}: \frac{D^2_y\widehat{ T_{\eps}v_{\eps}}}{\eps^{3/2}}d\bar xd y\eto  \int_{W\times Y\times (-\infty,0)}D^2_y\hat v(\bar x, y) : D^2_y\psi(\bar x, y) d\bar xdy\nonumber 
\end{eqnarray}
\hfill $\Box$\\

\begin{thm}\label{Lemma-micro3} {\bf (Characterization of $\hat v$ via  a two scale weak problem)} Let 
$\hat v\in  L^2(W, w^{2,2}_{Per_Y}(Y\times (-\infty , 0)))$ be the function from Theorem~\ref{mainmacroweak}. Then 
 \begin{equation}\label{micro5}
 \int_{W\times Y\times (-\infty,0)}D^2_y\hat v(\bar x, y): D^2_y\psi(\bar x, y)d\bar x dy=0\, ,
\end{equation}
for all  $\psi \in  L^2(W, w^{2,2}_{Per_Y}(Y\times (-\infty,0)))$ such that $\psi (\bar x, \bar y, 0)=0$ on $ W\times Y$. Moreover,
for any $i=1,\ldots, N-1$, we have 
\begin{equation}\label{Lemma-micro3bd}
\frac{\partial \hat v}{\partial y_i} (\bar x, \bar y,0)=\frac{\partial b}{\partial y_i}(\bar y)\frac{\partial v}{\partial x_N}(\bar x, 0), \ \ {\rm on}\ W\times Y\, .
\end{equation}
\end{thm}

 {\bf Proof.} The proof of (\ref{micro5}) for smooth test functions $\psi$ follows by passing to the limit in equation (\ref{micro1}) and combining (\ref{micro1})-(\ref{micro4}). The general case involving test functions $\psi \in L^2(W, w^{2,2}_{Per_Y}(Y\times (-\infty , 0)))$ follows by an approximation argument which we skip for brevity (we only mention that in order to preserve the boundary condition at $y_N=0$, one can first extend a given test function $\psi $ by setting $\psi (\bar x, \bar y ,-y_N )=-\psi (\bar x, \bar y ,y_N )$   and then using convolution). 

We now prove (\ref{Lemma-micro3bd}).   Let $V_{\eps}^{(i)}=\left(0, \dots, 0 , -\frac{\partial v_{\eps}}{\partial x_N},0, \dots , 0,   \frac{\partial v_{\eps}}{\partial x_i}\right)$ be the vector as in (\ref{casadovector}) for any $i=1, \dots , N-1$. As already noted, 
  $V_{\eps}^{(i)}\cdot \nu =0 $ on $\Gamma_{\eps}$ where $\nu$ is the normal to the boundary. 
We note that  by Lemma~\ref{unfoldconv} we have that 
$$
\frac{1}{\sqrt{\eps}}\frac{\partial }{\partial y_j}\left( \widehat{\frac{\partial v_{\eps} }{\partial x_i}}    \right)\weto 
\frac{\partial^2\hat v }{\partial y_j\partial y_i}
$$
in $L^2(W\times Y\times ]-\infty , 0[)$. Thus, it suffices to apply Casado Diaz et al \cite[Lemma 4.3]{casado} 
to conclude.
\hfill $\Box$\\

We plan now to describe function $\hat v$ in a more explicit way by separating the variables $\bar x $ and $y$ and providing a classical formulation of the microscopic problem (\ref{micro5}).   
To do so,  we shall also need to perform some calculations based on standard integration by parts on  domains of the type $Y\times ({d},0)$ with ${d}<0$. In essence, the computations are the same computations required 
to prove a known Green-type formula for the biharmonic operator.
For the sake of clarity, we state such formula and we provide a short proof in the classical setting. We recall that, given a smooth open set $U$ and a smooth vector field $F:\partial U\to {\mathbb{R}}^N$,  the tangential divergence of $F$ is defined by ${\rm div}_{\partial U}F={\rm div}F-(\nabla F\cdot \nu)\nu$ where it is meant that the vector field $F$ is extended smoothly in a neighborhood of $\partial U$. In the following statement, we shall also denote by  $F_{\partial U}$ the tangential  component of a vetor field $F$ as above, defined by $F_{\partial U}=F- (F\cdot \nu )\nu $.

\begin{lem}\label{bihagreen}{\bf (Biharmonic Green Formula)} Let $U$ be a bounded open set in ${\mathbb{R}}^N$ of class $C^{0,1}$. Let $f\in C^4(\bar U)$ and $\varphi \in C^2(\bar U)$. Then

\begin{eqnarray}\label{Green1Lip}
\lefteqn{\int_{U}D^2f: D^2\varphi dy  =
\int_{U}\Delta^2 f \varphi dy-\int_{\partial U}\frac{\partial \Delta f}{\partial \nu}\varphi  d\sigma +\int_{\partial U} \frac{\partial^2 f}{\partial \nu^2} \frac{\partial \varphi }{\partial \nu}   d\sigma } \nonumber \\
& &\qquad\qquad\qquad\qquad\qquad\qquad\qquad\qquad  + \int_{\partial U}(D^2 f \cdot \nu )_{\partial U}\cdot \nabla_{\partial U} \varphi d\sigma .
\end{eqnarray}

Moreover, if $\Omega $ is also of class $C^2$ then 
\begin{equation}\label{Green}
\int_{U}D^2f: D^2\varphi dy=  \int_U(\Delta^2 f)\varphi dy+\int_{\partial U}\frac{\partial^2 f}{\partial \nu^2}\frac{\partial \varphi }{\partial \nu}d\sigma -\int_{\partial U}
\left( {\rm div}_{\partial U}((D^2 f)\cdot \nu)_{\partial U}  +\frac{\partial \Delta f}{\partial \nu}\right)\varphi  d\sigma .
\end{equation}
 \end{lem}

{\bf Proof.}  Standard integration by parts  gives
\begin{eqnarray}\label{Green1}
\int_{U}D^2f: D^2 \varphi dy & =& \int_U \frac{\partial^2 f}{\partial y_i\partial y_j}\frac{\partial^2 \varphi }{\partial y_i\partial y_j}dy=-\int_U \frac{\partial^3 f}{\partial y_i^2\partial y_j}\frac{\partial \varphi }{\partial y_j}dy+\int_{\partial U} \frac{\partial^2 f}{\partial y_i\partial y_j}\frac{\partial \varphi }{\partial y_j}\nu_idy\nonumber\\
& = & 
\int_U \Delta f\Delta \varphi  dy-\int_{\partial U}\Delta f \frac{\partial \varphi }{\partial \nu }d\sigma 
+ \int_{\partial U}(D^2f \cdot \nu )\nabla \varphi  d\sigma\, , \nonumber \\
\end{eqnarray}
where  summation symbols have been dropped.
Formula (\ref{Green1Lip}) simply follows by the standard  Green formula and by decomposing  $\nabla \varphi$ as $\nabla_{\partial U} \varphi +\frac{\partial \varphi }{\partial \nu}\nu$ in (\ref{Green1}).  Finally, formula 
(\ref{Green}) follows by applying the Tangential Green Formula (see, e.g., Delfour and Zolesio~\cite[\S 5.5]{delzol})  to the last integral in the right-hand side of (\ref{Green1Lip}). \hfill $\Box$\\

We also need the following lemma where by $W^{4,2}_{Per_Y}(Y\times ({d},0))$  we denote the space  of functions in $W^{4,2}_{loc}({\mathbb{R}}^N\times (d , 0) )$  $Y$-periodic in the first $(N-1)$-variables $\bar y$ and such that all derivatives up to the fourth order are square summable in $Y\times ({d},0)$.

\begin{lem}\label{microy} There exists  $V\in   w^{2,2}_{Per_Y}(Y\times (-\infty , 0))$ satisfying the equation
 \begin{equation}\label{microweak}
 \int_{ Y\times (-\infty,0)}D^2V: D^2\psi dy=0  ,
\end{equation}
for all $\psi \in  w^{2,2}_{Per_Y}(Y\times (-\infty , 0))$ with $\psi ( \bar y, 0)=0$ on  $  Y$, and the boundary condition
\begin{equation}\label{microweak1} V (\bar y,0)=b(\bar y), \ \ {\rm on}\ Y.
\end{equation} 
Function $V$ is unique up to the sum of a monomial of the type $ay_n$ where  $a$ is any real number.  


Moreover,  $V$ is of class $W^{4,2}_{Per_Y}(Y\times ({d},0))$ for any ${d}<0$, 
$\Delta^2V= 0$  and it satisfies the  boundary condition
\begin{equation}\label{microweak2}
\frac{\partial ^2V}{\partial \nu^2}(\bar y, 0)=0,\ \ {\rm on}\ Y.
\end{equation}
\end{lem}

{\bf Proof. } In order to prove the existence of $V$, we use a standard direct method in the Calculus of Variations. 
Let $A=\{ V\in w^{2,2}_{Per_Y}(Y\times (-\infty ,0)):\ V (\bar y ,0)= b(\bar y)\ {\rm on }\ Y\}$. Obviously, $A\ne \emptyset$. It is clear that a minimizer (if it exists) for the problem 
\begin{equation}\label{min}
\inf_{V\in A}\int_{Y\times (-\infty ,0) }|D^2V|^2dy
\end{equation}
is a solution to problem (\ref{microweak}) and satisfies the required boundary condition (see also Lemma~\ref{bihagreen}). Thus it suffices to prove that there exists a minimizer in (\ref{min}). 
Let $V_n\in A$, $n\in {\mathbb{N}}$ be a minimizing sequence for (\ref{min}). Since this sequence is bounded in $w^{2,2}_{Per_Y}(Y\times (-\infty ,0))$ there exists 
$\tilde V\in w^{2,2}_{Per_Y}(Y\times (-\infty ,0))$ such that $D^2V_n$ converges weakly to $D^2\tilde V $ in $L^{2}(Y\times (-\infty ,0))$. Since the trace operator is compact we have that $V_n(\bar y, 0)$ converges strongly to $\tilde V (\bar y, 0)$ in $L^2(Y)$, hence  $\tilde V (\bar y, 0)=b(\bar y) $ on $Y$
and $V \in A$. We now prove that $\tilde V$ is in fact a minimizer. By the elementary inequality $|\alpha|^2\geq |\beta|^2+2\beta\cdot (\alpha-\beta)$ valid for all vectors $\alpha, \beta$ in any Euclidean space, we get
$$
\int_{Y\times (-\infty ,0) }|D^2V_ n|^2dy\geq \int_{Y\times (-\infty ,0) }|D^2 \tilde V|^2dy+ 2 \int_{Y\times (-\infty ,0) }D^2\tilde V: (D^2V_n-D^2\tilde V)dy.
$$
Passing to the limit in the previous inequality, using the weak convergence and the fact that the sequence is minimizing in (\ref{min}), we find out that 
$$
\inf_{u\in A}\int_{Y\times (-\infty ,0) }|D^2V|^2dy= \int_{Y\times (-\infty ,0) }|D^2\tilde V|^2dy,
$$
hence $\tilde V$ is a minimizer. 

Uniqueness is easily proved by observing that if  $V\in   w^{2,2}_{Per_Y}(Y\times (-\infty , 0))$ satisfies equation (\ref{microweak}) and the boundary condition $V(\bar y, 0)=0$ on $Y$, then $V$ itself can be tested in  (\ref{microweak}), and it easily follows that $D^2V= 0$. This implies that $V$ is a polynomial of the first degree of the type $\sum_{j=1}^Na_jy_j$. By the periodicity in $\bar y$ we get $a_j=0$ for all $j=1,\dots , N-1$ and the proof of the first part of the statement is complete. 

Regularity of $V$ is standard.   The rest of the proof follows by using in the weak formulation (\ref{microweak}) test functions $\psi $ as in the statement, with bounded support in the $y_N$-direction and using formula (\ref{Green1Lip}).  In fact, using such test functions $\psi$, we get  that the boundary terms corresponding to the subset 
$\partial Y\times (-\infty , 0)$     of the boundary cancel out because of  periodicity, hence 
$$
\int_{ Y\times (-\infty,0)}D^2V: D^2\psi dy = \int_{ Y\times (-\infty,0)}(\Delta^2 V)\psi dy+\int_{Y}\frac{\partial^2 V}{\partial \nu^2}\frac{\partial \psi}{\partial \nu}d\sigma =0.
$$
By the previous formula and the arbitrary choice of $\psi$ we deduce that $V$ is biharmonic and satisfies the boundary condition (\ref{microweak2}). \hfill $\Box$\\

\begin{lem}\label{finallem} Let $V$ be as in Lemma~\ref{microy}.  Then 
\begin{equation}\label{finalcor0}
\int_{Y\times (-1,0)}D^2V: D^2 (b(\bar y)(1+y_N)^3)dy=\int_{W\times ]-\infty , 0[}|D^2V|^2dy=-\int_Yb\frac{\partial }{\partial y_N}(\Delta_{\bar y} V+\Delta V)d\bar y.
\end{equation}
\end{lem}

{\bf Proof. } Let $\tilde\psi $ be the real-valued function defined on $Y\times ]-\infty , 0[$ by
$$
\tilde \psi (y)=\left\{\begin{array}{ll}b(\bar y)(1+y_N)^3,\ & \ {\rm if}\ (\bar y, y_N)\in Y\times [-1,0),\\
0,\ & \ {\rm if}\ (\bar y, y_N)\in Y\times ]-\infty ,-1).
\end{array}\right.
$$
Clearly, $\tilde \psi \in w^{2,2}_{Per_Y}(Y\times (-\infty,0)))$ and $\tilde \psi (\bar y, 0)= V(\bar y, 0)$. The proof of the first equality in (\ref{finalcor0}) then easily follows by using $\psi = V-\tilde \psi  $ 
as a test function in (\ref{microweak}). 

In order to prove  the second equality in (\ref{finalcor0}), we use formula  (\ref{Green1Lip}) with $U=Y\times (-1,0)$, $f=V$ and $\varphi (\bar y, y_N)= b(\bar y)(1+y_N)^3$. First of all, we note that 
$$
 ((D^2f)\cdot \nu)_{\partial U}  =\nabla_{\bar y}\left(\frac{\partial V}{\partial y_N} \right)
$$
on the part of the boundary given by $Y\times \{0\}$.  
 Then  by  Lemma~\ref{microy} and formula (\ref{Green1Lip}), by   exploiting the periodicity of $V$ and $b$  which allows to get rid of the boundary terms 
corresponding to the subset of the boundary $\partial Y\times (-1,0)$, we have that  
$$\int_{Y\times (-1,0)}D^2V: D^2 (b(\bar y)(1+y_N)^3)dy = \int_Y \nabla_{\bar y}\left(\frac{\partial V (\bar y, 0)}{\partial y_N} \right)\cdot  \nabla_{\bar y} b(\bar y) d\bar y
-\int_{Y}b\frac{\partial \Delta v}{\partial y_N}d\bar y
. $$
Integrating by parts the right-hand side of  the previous equality and using again the periodicity of the functions we conclude.   \hfill $\Box $\\

Finally, we can prove the following

\begin{thm}\label{veryfinal}
Let $V$ be as in Lemma~\ref{microy}. Let $v, \hat v$ be the functions defined in Theorem~\ref{mainmacroweak}. Then 
$\hat v (\bar x, y)= V(y)\frac{\partial v}{\partial x_N}(\bar x, 0) +a(\bar x)y_N$ for almost all $(\bar x, y)\in W\times Y\times ]-\infty , 0]$ where   $a(\cdot)\in L^2(W)$ is a  function  depending only on $\bar x$.  

Moreover, for the strange term in (\ref{mainmacroweak1}) we have 
\begin{eqnarray}\label{veryfinal1}\lefteqn{
\int_{W} \int_{ Y\times (-1,0)} D^2_y \hat v(\bar x, y) : D^2_{y}( b(\bar y)(1+y_N)^3  ) dy   \frac{\partial \varphi}{\partial x_N}(\bar x, 0)d\bar x}\nonumber \\
&=& \int_{W\times ]-\infty , 0[}|D^2V|^2dy\int_W\frac{\partial v}{\partial x_N}(\bar x, 0)d\bar x   \frac{\partial \varphi}{\partial x_N}(\bar x, 0)d\bar x \nonumber \\
&=& \gamma  \int_W\frac{\partial v}{\partial x_N}(\bar x, 0)d\bar x   \frac{\partial \varphi}{\partial x_N}(\bar x, 0)d\bar x,
\end{eqnarray}
where $\gamma$ is defined in (\ref{maingamma}). 
\end{thm}

{\bf Proof.}
Recall that function $\hat v$ satisfies problem (\ref{micro5}) and the boundary condition (\ref{Lemma-micro3bd}). By proceeding exactly as in the proof of Lemma~\ref{microy} one can easily see that such function $\hat v$ is unique up to the sum  of a function of the type $a(\bar x)y_N$ as in the statement. The proof then follows simply by observing that a function of the type $ V(y)\frac{\partial v}{\partial x_N}(\bar x, 0) $  as in the statement  satisfies   problem (\ref{micro5}) and the boundary condition (\ref{Lemma-micro3bd}). 
Equality (\ref{veryfinal1}) follows by (\ref{finalcor0}). 
\hfill $\Box$ \\

{\bf Acknowledgments.} 
The authors wish to express their gratitude to V.G. Maz'ya, V.I. Burenkov and J. Casado-D\'iaz for the always helpful discussions on this topic. 

The first author is partially supported by grant MTM2012-31298 from MINECO, Spain and from ``Programa de Financiaci\'on Grupos de Investigaci\'on UCM-Banco de Santander'', Grupo de Investigaci\'on CADEDIF- 920894.  

The second author acknowledges financial support from 
the research project `Singular perturbation problems for differential operators' Progetto di Ateneo of the University of Padova.
The second author is also member of the Gruppo Nazionale per l'Analisi Matematica, la Probabilit\`{a} e le loro Applicazioni (GNAMPA) of the
Istituto Nazionale di Alta Matematica (INdAM).

Both authors acknowledge the warm hospitality received by each other institution on the  occasion of several research visits.


\noindent Jos\'{e} M. Arrieta\\
Departamento de Matem\'{a}tica Aplicada\\
Universidad Complutense de Madrid\\
28040 Madrid, Spain\\
and \\
Instituto de Ciencias Matem\'aticas CSIC-UAM-UC3M-UCM \\
C/Nicolas Cabrera 13-15 \\ 
Cantoblanco 28049, Madrid. Spain. \\
e-mail: arrieta@mat.ucm.es\\

\noindent Pier Domenico Lamberti\\
Dipartimento di Matematica\\
Universit\`{a} degli Studi di Padova\\
Via Trieste 63\\ 
35121 Padova, Italy\\
e-mail: lamberti@math.unipd.it


\begin{thebibliography}{10}


\bibitem{arrJDE95} J.M. Arrieta, Neumann eigenvalue problems on exterior perturbation of the domain,  Journal of Differential Equations 118  (1995), no. 1, 54-103.


\bibitem{arrbrus2010} J.M. Arrieta, S.M Bruschi, 
Very rapidly varying boundaries in equations with nonlinear boundary conditions. The case of non uniform Lispschitz deformation, 
{\sl Discrete and Continuous Dynamical Systems B},  Volume 14, Number 2, pp. 327-351 (2010)  


\bibitem{arca} J.M. Arrieta and A.N. Carvalho,  Spectral convergence and nonlinear dynamics of reaction-diffusion equations under perturbations of the domain. J. Differential Equations 199, 2004, no. 1, 143-178.

\bibitem{ArrCarLoz} J.M. Arrieta, A.N Carvalho, G. Lozada-Cruz , Dynamics in
dumbbell domains I. Continiuity of the set of equilibria, J. Differential Equations 231 (2006),  551-597.

\bibitem{ArrLam} J.M. Arrieta, P.D. Lamberti, Spectral stability results for higher order operators under perturbations of the domain, C. R. Acad.Sci.Paris,  Ser.I 351(2013), 725-730


\bibitem{bucbut} D. Bucur and G. Buttazzo, Variational methods in some shape optimization problems. Appunti dei Corsi Tenuti da Docenti della Scuola. [Notes of Courses Given by Teachers at the School] Scuola Normale Superiore, Pisa, 2002. 217 pp.

\bibitem{Babuska}  I. Babu\v{s}ka, The theory of small changes in the domain of existence in the theory of Partial Differential Equations and its applications, 1963 Differential Equations and Their Applications (Proc. Conf., Prague, 1962) pp. 13-26 Publ. House Czechoslovak Acad. Sci., Prague; Academic Press, New York.

\bibitem{babvyb} I. Babu\v{s}ka, R. V\'{y}born\'{y}, Rudolf Continuous dependence of eigenvalues on the domain. Czechoslovak Math. J. 15 (90) 1965, 169-178. 

\bibitem{buoso} D. Buoso, P.D. Lamberti, Eigenvalues of polyharmonic operators on variable domains,  ESAIM Control Optim. Calc. Var. 19 (2013), no. 4, 1225-1235.

\bibitem{buosoplates} D. Buoso, P.D. Lamberti, Shape deformation for vibrating hinged plates,  Math. Methods Appl. Sci. 37 (2014), no. 2, 237-244.

\bibitem{bu}
V.I.~Burenkov, {Sobolev spaces on domains}, B.G. Teubner, Stuttgart, 1998.


\bibitem{bula2012} V. Burenkov, P.D. Lamberti, Sharp spectral stability estimates via the Lebesgue
measure of domains for higher order elliptic operators,   Rev. Mat. Complut. 25  (2012), 435-457.

\bibitem{bulahigh} V. Burenkov, P.D. Lamberti, Spectral stability of higher order uniformly
elliptic operators. Sobolev spaces in mathematics. II, 69-102, Int. Math. Ser.
(N. Y.), 9, Springer, New York, 2009.


 \bibitem{bulala} V. Burenkov, P.D. Lamberti and M. Lanza de Cristoforis, Spectral stability of nonnegative selfadjoint operators,  Sovrem. Mat. Fundam. Napravl.  15 (2006),  76-111  (in Russian. English transl. in J. Math. Sci. (N.Y.)  149 (2008), 1417-1452).

\bibitem{CarPis}A.N. Carvalho, S. Piskarev, A general approximation
scheme for attractors of abstract parabolic problems,
Numer. Funct. Anal. Optim.  27 (2006), 785-829.

\bibitem{casado}J. Casado-D\'iaz M. Luna-Laynez, F.J. Su\'arez-Grau,   Asymptotic Behavior of the Navier-Stokes System in a Thin Domain with Navier Condition on a Slightly Rough Boundary, SIAM J. Math. Anal. Vol 45, 3 (2013), 1641-1674.

\bibitem{casado10} Casado-D\'iaz, J.; Luna-Laynez, M.; Su\'arez-Grau, F. J. Asymptotic behavior of a viscous fluid with slip boundary conditions on a slightly rough wall. Math. Models Methods Appl. Sci. 20 (2010), no. 1, 121-156. 

\bibitem{chas} L.M. Chasman,  
 An isoperimetric inequality for fundamental tones of free plates. Comm. Math. Phys. 303 (2011),  no. 2, 421-449.

\bibitem{ciodamgri}  D. Cioranescu, A. Damlamian, G. Griso, The periodic unfolding method in homogenization. SIAM J. Math. Anal.  40  (2008),  no. 4, 1585-1620. 


\bibitem{ciora-murat} D. Cioranescu, F. Murat, 
Un terme \'{e}trange venu d'ailleurs. {\sl Nonlinear partial differential equations and their applications. Coll\`ege de France Seminar}, Vol. II (Paris, 1979/1980), pp. 98-138, 389-390, Res. Notes in Math., 60, Pitman, Boston, Mass.-London, 1982.

\bibitem{dan} D. Daners, Domain perturbation for linear and semi-linear boundary value problems. Handbook of differential equations: stationary partial differential equations. Vol. VI, 1-81, Handb. Differ. Equ., Elsevier/North-Holland, Amsterdam, 2008.


\bibitem{da} E.B.~Davies, {  Spectral theory and differential operators},
Cambridge University Press, Cambridge, 1995.

\bibitem{delzol} Delfour~M.C., Zolesio~J.-P., Shapes and geometries. Analysis, differential calculus, and optimization.
 Advances in Design and Control, 4. Society for Industrial and Applied Mathematics (SIAM), Philadelphia, PA, 2001. 

\bibitem{gaz} F. Gazzola, H-C. Grunau and G. Sweers, Polyharmonic boundary value problems. Positivity preserving and nonlinear higher order elliptic equations in bounded domains. Lecture Notes in Mathematics, 1991. Springer-Verlag, Berlin, 2010.

\bibitem{hale} J.K. Hale, Eigenvalues and perturbed domains. Ten mathematical essays on approximation in analysis and topology, 95-123, Elsevier B. V., Amsterdam, 2005.

\bibitem{henry} D. Henry, Perturbation of the boundary in boundary-value problems of partial differential equations. With editorial assistance from Jack Hale and Antonio Luiz Pereira. London Mathematical Society Lecture Note Series, 318. Cambridge University Press, Cambridge, 2005.

\bibitem{henrot} A, Henrot, Extremum problems for eigenvalues of elliptic operators. Frontiers in Mathematics. Birkh\"{a}user Verlag, Basel, 2006.


\bibitem{komo}  C. Komo, Influence of surface roughness to solutions of the Boussinesq equations with Robin boundary condition,  {\sl Rev Mat Complut.} 28, (2015), 123-155.

\bibitem{Maz-Naz} V.G. Maz'ya, S.A. Nazarov,  Paradoxes of limit passage in solutions of boundary value problems involving the approximation of smooth domains by polygonal domains,  Math. USSR Izvestiya 29 (1987), no. 3, 511-533.

\bibitem{Maz-Naz-Pla} V.G. Maz'ya, S.A. Nazarov,  B. Plamenevskii,  Asymptotic Theory of Elliptic
Boundary Value Problems in Singularly Perturbed Domains I and II,  Birkh\"auser,
Basel, 2000.

\bibitem{nec} J. Nec\v{a}s, 
 Les m\'{e}thodes directes en th\'{e}orie des \'{e}quations elliptiques. Masson et Cie, \'{E}diteurs, Paris, 1967.

\bibitem{nire} L. Nirenberg,  On elliptic partial differential equations. Ann. Scuola Norm. Sup. Pisa (3) 13 1959, 115-162.

\bibitem{Stu} F. Stummel,  Perturbation of domains in elliptic boundary-value problems, in: Lecture Notes Math., Vol. 503, Springer-Verlag, Berlin-Heidelberg-New York (1976), pp. 110-136.

 \bibitem{Vai77a} G. Vainikko,   \"Uber die Konvergenz und Divergenz von N\"aherungsmethoden bei Eigenwertproblemen, Math. Nachr. 78 (1977), 145-164.
 
\bibitem{Vai79} G. Vainikko, Regular convergence of operators and the approximate solution of equations, Mathematical Analysis, Vol. 16 (Russian),
5-53, 151, VINITI, Moscow, 1979.

\bibitem{Stein} E.M. Stein, Singular Integrals and Differentiability Properties of Functions, {\sl Princeton University Press}, Princeton, New Jersey, (1970).


\end{thebibliography}
\end{document}